\documentclass[12pt]{article}
\usepackage{fullpage}
\usepackage[english]{babel}
\usepackage{amsmath,amsthm,amssymb}
\usepackage{amsfonts}
\usepackage{amsbsy}
\usepackage{algorithm}
\usepackage{graphicx}
\usepackage{dsfont}
\usepackage{color}
\usepackage{mathrsfs}
\usepackage{bm}
\usepackage[flushleft]{threeparttable}
\usepackage{relsize}
\usepackage{natbib}
 \bibpunct[, ]{(}{)}{,}{a}{}{,}%
\newtheorem{theorem}{Theorem}[section]

\newtheorem{lemma}{Lemma}
\newtheorem{proposition}{Proposition}
\newtheorem{definition}{Definition}

\newtheorem{assumption}{Assumption}
\numberwithin{theorem}{section}
\numberwithin{corollary}{section}
\numberwithin{lemma}{section}
\numberwithin{proposition}{section}
\numberwithin{definition}{section}
\numberwithin{remark}{section}
\numberwithin{equation}{section}
\numberwithin{table}{section}

\newcommand{\enorm}[1]{\left\lVert#1\right\rVert_2}
\newcommand{\abs}[1]{\left|#1\right|}

\allowdisplaybreaks

\makeatletter
\let\save@mathaccent\mathaccent
\newcommand*\if@single[3]{%
  \setbox0\hbox{${\mathaccent"0362{#1}}^H$}%
  \setbox2\hbox{${\mathaccent"0362{\kern0pt#1}}^H$}%
  \ifdim\ht0=\ht2 #3\else #2\fi
  }
\newcommand*\rel@kern[1]{\kern#1\dimexpr\macc@kerna}
\newcommand*\widebar[1]{\@ifnextchar^{{\wide@bar{#1}{0}}}
{\wide@bar{#1}{1}}}
\newcommand*\wide@bar[2]{\if@single{#1}{\wide@bar@{#1}{#2}{1}}
{\wide@bar@{#1}{#2}{2}}}
\newcommand*\wide@bar@[3]{%
  \begingroup
  \def\mathaccent##1##2{%
    \let\mathaccent\save@mathaccent
    \if#32 \let\macc@nucleus\first@char \fi
    \setbox\z@\hbox{$\macc@style{\macc@nucleus}_{}$}%
    \setbox\tw@\hbox{$\macc@style{\macc@nucleus}{}_{}$}%
    \dimen@\wd\tw@
    \advance\dimen@-\wd\z@
    \divide\dimen@ 3
    \@tempdima\wd\tw@
    \advance\@tempdima-\scriptspace
    \divide\@tempdima 10
    \advance\dimen@-\@tempdima
    \ifdim\dimen@>\z@ \dimen@0pt\fi
    \rel@kern{0.6}\kern-\dimen@
    \if#31
      \overline{\rel@kern{-0.6}\kern\dimen@\macc@nucleus\rel@kern{0.4}
      \kern\dimen@}%
      \advance\dimen@0.4\dimexpr\macc@kerna
      \let\final@kern#2%
      \ifdim\dimen@<\z@ \let\final@kern1\fi
      \if\final@kern1 \kern-\dimen@\fi
    \else
      \overline{\rel@kern{-0.6}\kern\dimen@#1}%
    \fi
  }%
  \macc@depth\@ne
  \let\math@bgroup\@empty \let\math@egroup\macc@set@skewchar
  \mathsurround\z@ \frozen@everymath{\mathgroup\macc@group\relax}%
  \macc@set@skewchar\relax
  \let\mathaccentV\macc@nested@a
  \if#31
    \macc@nested@a\relax111{#1}%
  \else
    \def\gobble@till@marker##1\endmarker{}%
    \futurelet\first@char\gobble@till@marker#1\endmarker
    \ifcat\noexpand\first@char A\else
      \def\first@char{}%
    \fi
    \macc@nested@a\relax111{\first@char}%
  \fi
  \endgroup
}
\makeatother

\begin{document}

\title{Simulation Optimization of Risk Measures with Adaptive Risk Levels}
\author{Helin Zhu, Joshua Hale, and Enlu Zhou\\
H. Milton Stewart School of Industrial and Systems Engineering,\\
Georgia Institute of Technology}
\date{}
\maketitle

\begin{abstract}
Optimizing risk measures such as Value-at-Risk (VaR) and Conditional Value-at-Risk (CVaR) of a general loss distribution is usually difficult, because
1) the loss function might lack structural properties such as convexity or differentiability since it is often generated via black-box simulation of a stochastic system;
2) evaluation of risk measures often requires rare-event simulation, which is computationally expensive.
In this paper, we study the extension of the recently proposed gradient-based adaptive stochastic search (GASS) to the optimization of risk measures VaR and CVaR. Instead of optimizing VaR or CVaR at the target risk level directly, we incorporate an adaptive updating scheme on the risk level, by initializing the algorithm at a small risk level and adaptively increasing it until the target risk level is achieved while the algorithm converges at the same time.
This enables us to adaptively reduce the number of samples required to estimate the risk measure at each iteration, and thus improving the overall efficiency of the algorithm.
\end{abstract}

\section{Introduction}\label{sec1:Introduction}

Risk measures such as Value-at-Risk (VaR) and Conditional Value-at-Risk (CVaR) are widely studied in various fields, in order to quantify the extreme behaviors of the distributions of interest.  Loosely speaking, VaR characterizes a tail quantile of a distribution, and CVaR characterizes the conditional expectation of the tail portion of the distribution.
VaR, as one of the earliest risk measures introduced in financial risk management, is easy to understand and interpret for practitioners. CVaR, as a classic coherent risk measure (see, e.g., \cite{artzner1997coherent}), exhibits nice properties such as convexity and monotonicity for optimization.
An abundant literature has dedicated to studying the estimation and optimization of risk measures under various settings. \cite{rockafellar2000optimization} and
\cite{rockafellar2002conditional} derive some fundamental properties of CVaR for general loss distributions in finance, and propose the fundamental minimization formula to facilitate the optimization of CVaR.
\cite{ruszczynski2006optimization} develop a dual representation for optimization of general coherent risk measures, and derive the optimality conditions via the dual representation. \cite{ruszczynski2010risk} study the optimization of risk measures under a multistage setting, and propose a risk-averse dynamic programming approach to risk optimization in Markov decision processes.
\cite{alexander2006minimizing} study the optimization of VaR and CVaR for derivatives portfolios with the addition of a cost that is proportional to the portfolio position.

In general, optimizing risk measures over continuous decision variables is a challenging problem, especially when the underlying loss function does not possess good structural properties such as convexity or differentiability.
Traditional gradient-based optimization methods often are not applicable, since little problem-specific knowledge is available when the loss function is evaluated via black-box simulation of a stochastic system.
In contrast, model-based optimization methods are good alternatives as they impose minimal requirements on problem structure. Examples of model-based methods include but are not limited to ant colony optimization (\cite{dorigo2005ant}),
annealing adaptive search (AAS) (\cite{romeijn1994simulated}),
the estimation of distribution algorithms (EDA) (\cite{larranaga2002estimation}),
the cross-entropy (CE) method (\cite{rubinstein2001combinatorial}), model reference adaptive search (MRAS) (\cite{hu2007model} and \cite{hu2008model}),
the interacting-particle algorithm (\cite{molvalioglu2009interacting}, \cite{molvalioglu2010meta}), and
gradient-based adaptive stochastic search (GASS) (\cite{zhou2014gradient}).

The main idea of model-based methods is to introduce a sampling distribution, which often belongs to a parameterized family of densities, over the solution space, and iteratively update the sampling distribution (or its parameter) by drawing and evaluating candidate solutions according to the sampling distribution.
The hope is to have the sampling distribution more and more concentrated on the promising region of the solution space where the optimal solutions are located, and eventually become a degenerate distribution on one of the global optima.
Therefore, finding an optimal solution in the solution space is transformed to finding an optimal sampling distribution parameter in the parameter space.
A key difference among the aforementioned model-based methods lies in how to update the sampling distribution.
For example, in CE and MRAS, the updating rule is derived by minimizing the Kullback-Leibler (K-L) divergences between a converging sequence of reference distributions and a chosen exponential family of densities.
For another example, in GASS, the updating rule on the sampling distribution parameter is derived by converting the original (possibly non-differentiable) deterministic optimization problem into a differentiable stochastic optimization problem on the sampling distribution parameter, and then applying a Newton-like scheme.

Compared with traditional gradient-based methods, model-based methods are more robust in the sense that at every iteration they exploit the promising region of the solution space that has already been identified, while maintaining the exploration of the entire solution space.
The updating rule on the sampling distribution parameter controls the balance between exploration and exploitation.

Although all the aforementioned model-based methods are designed for deterministic optimization problems, they can be extended to risk (VaR or CVaR) optimization problems in which the exact risk values are replaced with (biased) sample estimates.
However, a straightforward extension usually leads to a computationally expensive algorithm, due to the rare-event simulation required in estimating the risk values.
This issue is even more severe when the risk level of interest is close to $1$, which is often the case for risk management practitioners.
It inspires us to consider the following question. For a risk optimization problem, is it possible to initialize a model-based algorithm at a small risk level (close to $0$), and then adaptively increase the risk level at every iteration such that
the target risk level is achieved while the algorithm converges at the same time?
The motivation is that the algorithm will consume less simulation budget (since the risk level is small) during the ``warm-up'' phase of the algorithm, and solve problems that are close to the original one during the ``convergence'' phase of the algorithm, eventually leading to total budget saving.
The key to this question lies in finding a signal to link the updating rule on the risk level with the updating rule on the sampling distribution parameter, where the signal is capable of measuring empirically the algorithm's emphasis between the exploitation of a promising region and the exploration of the entire solution space.

In this paper, we will focus on the extension of a specific model-based methods, i.e., GASS by \cite{zhou2014gradient}, to the optimization of risk measures.
We choose  GASS because it could also be interpreted as a gradient-based scheme of a reformulated problem, in which a Newton-like updating rule is applied on the sampling distribution parameter, and thus the gradient (even the Hessian) in the updating rule of the sampling distribution parameter can be viewed as a signal that empirically measures the algorithm's emphasis between the exploitation of a promising region and the exploration of the entire solution space.
Therefore, we could adjust the risk level adaptively using the information contained in the gradient (e.g., its norm) at every iteration.
In particular, we will propose an updating rule that increases the risk level proportionally to the decrease in the norm of the gradient. We will show that incorporating such an updating scheme on the risk level in the algorithm guarantees that the target risk level is achieved at the same time when the algorithm converges. Furthermore, compared with vanilla GASS, we will show that the proposed algorithm achieves significant total budget savings.

To the best of our knowledge, this work is among the first to apply a model-based algorithm to risk optimization problems, and among the first to propose a risk optimization scheme with adaptive risk levels.
To ease the presentation, we will only focus on CVaR optimization, and the extension of the proposed algorithm to VaR optimization (and possibly other risk measures such as probability of large loss) is straightforward.

The rest of the paper is organized as follows.
In Section \ref{sec2:Formulation}, we will first describe the CVaR optimization problem. Then we extend the algorithm GASS, which is originally developed for deterministic non-differentiable optimization problems, to the CVaR optimization problem.
The detailed algorithms are presented in Section \ref{sec3:Algorithm}, in which Algorithm \ref{alg.3.1} (referred to as ``GASS-CVaR'') is a straightforward extension of GASS and Algorithm \ref{alg.3.2} (referred to as ``GASS-CVaR-ARL'') further incorporates an updating rule that adaptively adjusts the risk level.
Convergence analysis of both algorithms are presented in Section \ref{sec4:Convergence}.
In Section \ref{sec5:Numerical}, we illustrate the performance of the proposed algorithms by carrying out numerical tests on several benchmark loss functions.
We conclude the paper in Section \ref{sec6:Conclusion}.

\section{General Framework}\label{sec2:Formulation}


Consider a scalar loss function of the form $l(x, \xi_x)$, where $x\in \mathcal{X}\subseteq \mathbb{R}^{d_x}$ represents the decision variable, and $\xi_x$ represents the randomness in the loss function. The distribution of $l(x, \xi_x)$ may or may not depend on $x$.
The loss function $l(x, \xi_x)$ can be evaluated either directly or through simulation.
Furthermore, to ease the presentation, we assume $l(x, \xi_x)$ admits an almost everywhere (a.e.) positive and continuous probability density function (p.d.f.) $p(t;x)$, and thus a continuous and strictly increasing cumulative distribution function (c.d.f.) $P(t;x)$ for all $x\in \mathcal{X}$.
The objective is to minimize CVaR of $l(x;\xi_x)$ at a risk level of interest $\alpha^\ast$ ($0<\alpha^\ast<1$) with respect to (w.r.t.) $x\in \mathcal{X}$.
That is, to solve the following stochastic optimization problem:
\begin{equation}\label{eq.2.1}
\min_{x\in \mathcal{X}}\;C_{\alpha^\ast}(x)\overset{\triangle}
=CVaR_{\alpha^\ast}\left(l(x, \xi_x)\right),\quad \mbox{or equivalently,}\quad \max_{x\in \mathcal{X}}\;-C_{\alpha^\ast}(x),
\end{equation}
where $CVaR_{\alpha^\ast}\left(l(x, \xi_x)\right)$ is defined by
\begin{equation*}
\begin{split}
CVaR_{\alpha^\ast}\left(l(x, \xi_x)\right)&\overset{\triangle}=
\mathbb{E}_{\xi_x}\left[l(x, \xi_x)|l(x, \xi_x)\ge V_{\alpha^\ast}(x)
\right]\\
&{}=\frac{1}{1-\alpha^\ast}\mathbb{E}_{\xi_x}\left[l(x, \xi_x)\mathds{1}
\left\{l(x, \xi_x)\ge V_{\alpha^\ast}(x)\right\}\right]\\
&{}=V_{\alpha^\ast}
(x)+\frac{1}{1-\alpha^\ast}\mathbb{E}_{\xi_x}\left[\left(
l(x, \xi_x)-V_{\alpha^\ast}(x)\right)^+\right],
\end{split}
\end{equation*}
where $\mathds{1}\{A\}$ is $1$ if event $A$ is true and $0$ otherwise, $(u)^+=\max(u, 0)$, and $V_{\alpha^\ast}(x)$ is VaR of $l(x, \xi_x)$ at the target risk level $\alpha^\ast$, i.e.,
\begin{equation*}
V_{\alpha^\ast}(x)=VaR_{\alpha^\ast}\left(l(x, \xi_x)\right)\overset{\triangle}=\inf\{t: P(t; x)\ge \alpha^\ast\}=P^{-1}(\alpha^\ast; x).
\end{equation*}
Note that the inverse c.d.f. $P^{-1}(\alpha^\ast; x)$ exists because $P(t; x)$ is strictly increasing in $t$.
We also follow the standard assumption that $C_{\alpha^\ast}(x)$ is bounded from below and above on $\mathcal{X}$, i.e., $\exists C_{lb}>-\infty, C_{ub}<\infty$ s.t. $C_{lb}<C_{\alpha^\ast}(x)<C_{ub},\;\forall x\in \mathcal{X}$.

Problem (\ref{eq.2.1}) might be difficult to solve when $l(x;\xi_x)$ lacks structural properties such as convexity and differentiability. Thus, traditional gradient-based method might not be applicable.
In contrast, model-based methods are good alternatives as in general they impose few requirements on the structure of $l(x;\xi_x)$. Therefore, we will apply a model-based method to solve problem (\ref{eq.2.1}).
In principle, we could extend the algorithm GASS in \cite{zhou2014gradient}, which is reviewed in next section.

\subsection{Review of GASS}\label{sec2.subsec1:Review}

Similar to many other model-based methods, the main idea of GASS is to introduce a parameterized sampling distribution over the solution space, and update the sampling distribution iteratively towards the promising region of the solution space.
Let us illustrate the main idea in a general framework, where one aims to maximize a deterministic function $L(x)$ over $x\in \mathcal{X}$.

Introduce a parameterized family of densities
$\{f(x;\theta): \theta\in \Theta \subset \mathbb{R}^{d_\theta}\}$ as the sampling distribution, where $\theta$ represents the parameter that will be updated over iterations. Consider a simple reformulation as follows:
$$
H(\theta)\overset{\triangle}=\int L(x)f(x; \theta) dx.
$$
Then $H(\theta)\le L(x^{\ast})=L^{\ast}$, where $x^{\ast}$ denotes the optimal solution or one of the optima, and $L^{\ast}$ denotes the optimal function value on $\mathcal{X}$.
Note that the equality is achieved if and only if all the probability mass of $f(x;\theta)$ concentrates on a subset of the set of global optima.
Given the existence of such a $\theta$, one could solve the reformulated problem $\max_{\theta\in \Theta}H(\theta)$ instead of the original problem, since the optimal parameter will recover the optimal solution and the optimal function value.

An advantage of the reformulated problem over the original problem is that it is differentiable in $\theta$ under mild regularity conditions on $f(x;\theta)$, and the gradient is easy to derive as follows:
$$
\bigtriangledown_{\theta}H(\theta)=\bigtriangledown_{\theta}
\int L(x)f(x; \theta) dx=\int L(x)
\frac{\bigtriangledown_{\theta}f(x; \theta)}{f(x; \theta)}
f(x; \theta) dx=
\mathbb{E}_{f(\cdot;\theta)}\left[L(x)\bigtriangledown_{\theta}\ln f(x;\theta)\right].
$$
Note that an unbiased estimator of $\bigtriangledown_{\theta}H(\theta)$ could be obtained by drawing independent and identically distributed (i.i.d.) samples $x^i \sim f(x;\theta), i=1,...,N$, evaluating $L(x^i)\bigtriangledown_{\theta}\ln f(x^i;\theta)$, and taking the sample average of $\{L(x^i)\bigtriangledown_{\theta}\ln f(x^i;\theta)\}$.
Therefore, one could solve the reformulated problem via a (stochastic) gradient-based method. Specifically, the method iteratively carries out the following two steps:
\begin{itemize}
\item[1.] Generate candidate solutions according to the sampling distribution.
\item[2.] Based on the evaluation of the candidate solutions, update the parameter of the sampling distribution via gradient search.
\end{itemize}
Intuitively, it combines the relative fast convergence of gradient search with the robustness of model-based optimization in terms of maintaining a global exploration of the solution space.


Based on the above main idea, now let us review the full-blown GASS algorithm.
Introduce a shape function $S_\theta: \mathbb{R}\rightarrow \mathbb{R}^+$, where the subscript $\theta$ signifies the possible dependence of the shape function on the parameter $\theta$. The shape function $S_\theta$ satisfies the following conditions:
for every $\theta$, $S_\theta(y)$ is strictly increasing in $y$, and bounded below from zero and above for finite $y$;
moreover, for every fixed $y$, $S_{\theta}(y)$ is continuous in $\theta$.
The purpose of introducing $S_\theta$ is to make the objective function positive, and yet preserve the order of the solutions and in particular the optimal solution. Moreover, the dependence of $S_\theta$ on $\theta$
adds flexibility to the algorithm by giving user the freedom to choose a weighting scheme on the samples based on the function evaluations. For example, a good choice of $S_{\theta}(\cdot)$ is
\begin{equation}\label{eq.2.2}
S_\theta(L(x))=\frac{1}{1+\exp(-S_o(L(x)-\gamma_\theta))},
\end{equation}
where $S_o$ is a large positive constant, and $\gamma_\theta$ is the $(1-\rho)$-quantile
\begin{equation}\label{eq.2.3}
\gamma_\theta\overset{\triangle}=\sup_{r}\left\{
r: P_{f(\cdot;\theta)}\left\{x\in \mathcal{X}: L(x)\ge r\right\} \ge \rho\right\},
\end{equation}
where $P_{f(\cdot;\theta)}\{A\}$ denotes the probability of event $A$ w.r.t. $f(\cdot;\theta)$. Notice that $S_{\theta}(\cdot)$ could be viewed as a continuous approximation of the indicator function $\mathds{1}\{L(x)\ge \gamma_\theta\}$ that gives equal weights to the solutions with function values above $\gamma_\theta$ and eliminates the solutions with function values below $\gamma_\theta$.

For an arbitrary but fixed $\theta^\prime \in \Theta$, define
\begin{equation}\label{eq.2.4}
H(\theta;\theta^\prime)\overset{\triangle}=\int S_{\theta^\prime}(L(x)) f(x;\theta)dx,\quad \mbox{and}
\quad h(\theta;\theta^\prime) \overset{\triangle}=\ln H(\theta;\theta^\prime).
\end{equation}
By the condition on the shape function and the fact that $\ln(\cdot)$ is a strictly increasing function, the original problem can be transform to $\max_{\theta\in \Theta} h(\theta;\theta^\prime)$ for any fixed $\theta^\prime$. Following the main idea outlined before, \cite{zhou2014gradient} propose a stochastic search algorithm that iteratively carries out the following two steps:

\begin{itemize}
\item[1.] Generate candidate solutions from $f(x;\theta_k)$, where $\theta_k$ is the sampling distribution parameter obtained at iteration $k$.
\item[2.] Update $\theta_k\rightarrow\theta_{k+1}$ using a Newton-like iteration for $\max_{\theta}h(\theta;\theta_k)$, where
    $h(\theta;\theta_k)=\ln\{\int S_{\theta_k}(L(x))f(x;\theta)dx\}$.
\end{itemize}

Note that the second step requires us to compute the gradient and Hessian of $h(\theta;\theta_k)$ at $\theta=\theta_k$, which, as shown by \cite{zhou2014gradient}, have analytical expressions as expectations under certain probability measures. In particular, if the sampling distribution belongs to an exponential family of densities, which is defined in the following Definition \ref{def.2.1}, then these expressions can be further simplified.
\begin{definition}\label{def.2.1}
A family $\{f(x;\theta): \theta\in \Theta\}$ is an exponential family of densities if it satisfies
\begin{equation*}
f(x;\theta)=\exp\left\{\theta^T \Gamma(x)-\eta(\theta)\right\},
\end{equation*}
where $\Gamma(x)=[\Gamma_1(x),...,\Gamma_d(x)]^T$ is the vector of sufficient statistics,
$$\eta(\theta)=\ln\{\int\exp(\theta^T \Gamma(x))dx\}$$
is the normalization factor that ensures $f(x;\theta)$ to be a p.d.f.; $\Theta=\{\theta: \abs{\eta(\theta)}<\infty\}$ is the natural parameter space with a nonempty interior.
\end{definition}
Proposition \ref{prop.2.1} below provides the corresponding analytical expressions of the gradient and Hessian of $h(\theta;\theta_k)$ at $\theta=\theta_k$, when an exponential family of densities is used as the sampling distribution. We refer to \cite{zhou2014gradient} for the detailed derivations.

\begin{proposition}\label{prop.2.1}
If $\{f(x;\theta): \theta\in \Theta\}$ is an exponential family of densities, then the gradient and Hessian of $h(\theta;\theta_k)$ at  $\theta=\theta_k$ have closed-form expressions as follows.
\begin{equation}\label{eq.2.5}
\left\{\begin{array}{l}
\bigtriangledown_{\theta}h(\theta;\theta_k)\big|
_{\theta=\theta_k}=\mathbb{E}_{q(\cdot;\theta_k)}
\left[\Gamma(x)\right]-\mathbb{E}_{\theta_k}
\left[\Gamma(x)\right],\vspace{2mm}
\\
\bigtriangledown^2_{\theta}h(\theta;\theta_k)\big|
_{\theta=\theta_k}=Var_{q(\cdot;\theta_k)}
\left[\Gamma(x)\right]-Var_{\theta_k}
\left[\Gamma(x)\right],
\end{array}
\right.
\vspace{-2mm}
\end{equation}
where
$$
q(x;\theta_k)= \frac{ S_{\theta_k}(L(x)) f(x;\theta_k)}{\int S_{\theta_k}(L(x)) f(x;\theta_k)dx}
$$
is a ``re-weighted'' p.d.f.; $\mathbb{E}_{q(\cdot;\theta_k)}
\left[\cdot\right]$ and $Var_{q(\cdot;\theta_k)}
\left[\cdot\right]$ denote the expectation and variance w.r.t. $q(\cdot;\theta_k)$, respectively; $\mathbb{E}_{\theta_k}
\left[\cdot\right]$ and $Var_{\theta_k}
\left[\cdot\right]$ denote the expectation and variance w.r.t. $f(\cdot;\theta_k)$, respectively.
\end{proposition}

Note that the Hessian $\bigtriangledown^2_{\theta}h(\theta;\theta_k)\big|
_{\theta=\theta_k}$ might not be negative semi-definite. To ensure the parameter updating is along the ascent direction of $h(\theta;\theta_k)$ in a Newton-like scheme, one could approximate $\bigtriangledown^2_{\theta}h(\theta;\theta_k)\big|
_{\theta=\theta_k}$ by a negative-definite term $-(Var_{\theta_k}[T(x)]+\epsilon I)$, which is a slight perturbation of the second term in $\bigtriangledown^2_{\theta}h(\theta;\theta_k)\big|
_{\theta=\theta_k}$. Here $\epsilon$ is a small positive number and $I$ is the identity matrix of proper dimension. Then, a Newton-like updating scheme of $\theta$ for $\max_{\theta}h(\theta;\theta_k)$ is as follows.
\begin{eqnarray}\label{eq.2.6}
\theta_{k+1}&=&\Pi_{\Theta}\left\{\theta_k+\beta_k\left(
Var_{\theta_k}[\Gamma(x)]+\epsilon I\right)^{-1}\bigtriangledown_{\theta}h(\theta;\theta_k)\big|
_{\theta=\theta_k}\right\}\nonumber\\
&=&\Pi_{\Theta}\left\{\theta_k+\beta_k\left(
Var_{\theta_k}[\Gamma(x)]+\epsilon I\right)^{-1}\left(
\mathbb{E}_{q_k}
\left[\Gamma(x)\right]-\mathbb{E}_{\theta_k}
\left[\Gamma(x)\right]\right)\right\},
\end{eqnarray}
where $\beta_k$ is a positive step-size, $\mathbb{E}_{q_k}
\left[\cdot\right]$ denotes the expectation w.r.t. $q(\cdot;\theta_k)$, and $\Pi_{\Theta}\{\cdot\}$ denotes the projection operator that projects an iterate back onto the parameter space $\Theta$ by choosing the closest point in $\Theta$.

To have an implementable algorithm, the expectation and variance terms in (\ref{eq.2.6}) need to be evaluated or estimated. Notice that the expectation term $\mathbb{E}_{\theta_k}[\Gamma(x)]$ can be calculated analytically in most cases. For example, if the chosen exponential family of densities is the Gaussian family, then $\mathbb{E}_{\theta_k}[\Gamma(x)]$ reduces to the mean and second moment of a Gaussian distribution. The variance term $Var_{\theta_k}[\Gamma(x)]$ might not be directly available, but it could be estimated by the sample variance using the candidate solutions drawn from $f(\cdot;\theta_k)$. Specifically, suppose $N_k$ i.i.d. samples $\{x^i_k: i=1,...,N_k\}$ are drawn from $f(x;\theta_k)$, then
\begin{equation}\label{eq.2.7}
\widehat{Var}_{\theta_k}\left[\Gamma(x)\right]\overset{\triangle}=
\frac{1}{N_k-1}\sum_{i=1}^{N_k}\Gamma(x_k^i)\Gamma(x_k^i)^T
-\frac{1}{N^2_k-N_k}\left(\sum_{i=1}^{N_k}
\Gamma(x_k^i)\right)\left(\sum_{i=1}^{N_k}\Gamma(x_k^i)\right)^T
\end{equation}
is a sample estimate of $Var_{\theta_k}[\Gamma(x)]$. The remaining term $\mathbb{E}_{q_k}[\Gamma(x)]$ can be estimated using the principle of importance sampling with samples $\{x^i_k\}$, noting that
$$
\mathbb{E}_{q_k}
\left[\Gamma(x)\right]\propto \int S_{\theta_k}(L(x))\Gamma(x) f(x;\theta_k)dx.
$$
That is, the expectation $\mathbb{E}_{q_k}[\Gamma(x)]$ could be estimated  by
\begin{equation}\label{eq.2.8}
\widetilde{\mathbb{E}}_{q_k}\left[\Gamma(x)\right]\overset{\triangle}=
\sum_{i=1}^{N_k}w_k^i \Gamma(x^i_k),
\end{equation}
where $\{w^i_k: i=1,...,N_k\}$ are self-normalized weights given by
$$
w_k^i=\frac{S_{\theta_k}(L(x_k^i))}
{\sum_{j=1}^{N_k}S_{\theta_k}(L(x_k^j))},\quad i=1,...,N_k.
$$
When $S_{\theta_k}(\cdot)$ takes a form such as (\ref{eq.2.2}), it has to be estimated by samples as well since the $(1-\rho)$-quantile $\gamma_{\theta_k}$ defined in (\ref{eq.2.3}) needs to be estimated by a sample $(1-\rho)$-quantile. Denote the sample quantile by $\widehat{\gamma}_{\theta_k}$ and the resulted approximate shape function by $\widehat{S}_{\theta_k}(\cdot)$. Then, $\{w^i_k\}$ are approximated according to
$$
\widehat{w}_k^i=\frac{\widehat{S}_{\theta_k}(L(x_k^i))}
{\sum_{j=1}^{N_k}\widehat{S}_{\theta_k}(L(x_k^j))},\quad i=1,...,N_k,
$$
and thus $\mathbb{E}_{q_k}
\left[\Gamma(x)\right]$ is approximated by
\begin{equation}\label{eq.2.9}
\widehat{\mathbb{E}}_{q_k}\left[\Gamma(x)\right]\overset{\triangle}=
\sum_{i=1}^{N_k}\widehat{w}_k^i \Gamma(x^i_k).
\end{equation}
Eventually, the gradient $g_k:= \mathbb{E}_{q_k}
\left[\Gamma(x)\right]-\mathbb{E}_{\theta_k}
\left[\Gamma(x)\right]$ is approximated by
\begin{equation*}
\widehat{g}_k\overset{\triangle}
=\widehat{\mathbb{E}}_{q_k}\left[\Gamma(x)\right]-
\mathbb{E}_{\theta_k}
\left[\Gamma(x)\right].
\end{equation*}

\subsection{Extension of GASS to Optimization of CVaR}

When CVaR of the loss function $l(x,\xi_x)$, $C_{\alpha^\ast}(x)$, could be evaluated exactly for all $x\in \mathcal{X}$, we can directly apply the scheme described above to solve the CVaR minimization problem (\ref{eq.2.1}).
Since $l(x,\xi_x)$ is usually evaluated via simulation, generally its p.d.f. and c.d.f. are not available; thus, $C_{\alpha^\ast}(x)$ could not be evaluated analytically.
Nevertheless, it could be estimated via Monte Carlo simulation. In particular, suppose $M$ i.i.d. loss samples $\{l(x,\xi_x^1), l(x,\xi_x^2), ..., l(x,\xi_x^M)\}$ are simulated, and then sorted in ascending order as
$
l(x,\xi_x^{(1)})\le l(x,\xi_x^{(2)})\le  ... \le l(x,\xi_x^{(M)}),
$
which forms an empirical loss distribution. A natural estimator of $C_{\alpha^\ast}(x)$ is CVaR of the empirical loss distribution, which is defined as follows.
\begin{equation}\label{eq.2.10}
\widehat{C}_{\alpha^\ast}(x)\overset{\triangle}=
\widehat{V}_{\alpha^\ast}(x)+
\frac{1}{M(1-\alpha^\ast)}\sum_{m=1}^{M}\left(l(x, \xi_x^m)-\widehat{V}_{\alpha^\ast}(x)\right)^+,
\end{equation}
where
\begin{equation}\label{eq.2.11}
\widehat{V}_{\alpha^\ast}(x)\overset{\triangle}=l\left(x, \xi_x^{(\lceil\alpha^\ast M\rceil)}\right)
\end{equation}
is VaR of the empirical loss distribution that plays the role of VaR estimator, and $\lceil\alpha^\ast M\rceil$ is the smallest integer that is greater than or equal to $\alpha^\ast M$.

Although the estimator $\widehat{C}_{\alpha^\ast}(x)$ is biased, it is strongly consistent and asymptotic normally distributed under mild regularity assumptions (see, e.g., \cite{zhu2016risk}). In principle, we can use it as a replacement for ${C}_{\alpha^\ast}(x)$ and plug it into GASS algorithm.

\section{Algorithms: GASS-CVaR,  GASS-CVaR-ARL}\label{sec3:Algorithm}

Now let us formally present the following Algorithm \ref{alg.3.1}, which is referred to as GASS-CVaR, for simulation optimization of CVaR.

\begin{algorithm}
\caption{\textbf{G}radient-based \textbf{A}daptive \textbf{S}tochastic \textbf{S}earch for Optimization of \textbf{CVaR}}\label{alg.3.1}
1. \textbf{Initialization}: Choose an exponential family of densities $\{f(x;\theta): \theta\in \Theta\}$, and specify a small positive constant $\epsilon$, initial parameter $\theta_0$, sample size sequence $\{N_k\}$ that satisfies $N_k \rightarrow\infty$, simulation budget sequence $\{M_k\}$ that satisfies $M_k \rightarrow\infty$, and step size sequence $\{\beta_k\}$ that satisfies $\sum_{k=0}^\infty\beta_k=\infty, \sum_{k=0}^\infty\beta^2_k<\infty$. Set $k=0$.
\\
2.\textbf{Sampling}: Draw candidate solutions $\{x_k^i \overset{i.i.d.}\sim f(x;\theta_k): i=1, 2,..., N_k\}$. For each $x_k^i$, simulate i.i.d. loss scenarios $\{l(x_k^i, \xi_k^{i,j}): j=1,...,M_k\}$, and sort them in ascending order, denoted by
$$
l\left(x_k^i, \xi_k^{i,(1)}\right)\le l\left(x_k^i, \xi_k^{i,(2)}\right)\le \cdots\le l\left(x_k^i, \xi_k^{i,(M_k)}\right).
$$
Estimate the CVaR of the loss for each candidate solution at target risk level $\alpha^\ast$:
$$
\widehat{C}_{\alpha^\ast}(x_k^i)=l\left(x_k^i, \xi_k^{i,(\lceil \alpha^\ast M_k\rceil)}\right)+\frac{1}{M_k(1-\alpha^\ast)}
\sum_{j=1}^{M_k}\left(l\left(x_k^i, \xi_k^{i,j}\right)-
l\left(x_k^i, \xi_k^{i,(\lceil \alpha^\ast M_k\rceil)}\right)\right)^+.
$$
3. \textbf{Estimation}: Compute the normalized weights $\widehat{w}_k^i$ according to
$$
\widehat{w}_k^i=\frac{\widehat{S}_{\theta_k}
\left(-\widehat{C}_{\alpha^\ast}(x_k^i)\right)}
{\sum_{j=1}^{N_k}\widehat{S}_{\theta_k}
\left(-\widehat{C}_{\alpha^\ast}(x_k^j)\right)}, \quad i=1,..., N_k,
$$
and then estimate $\mathbb{E}_{q_k}[\Gamma(x)]$ and $Var_{\theta_k}[\Gamma(x)]$ via
\begin{equation*}
\left\{\begin{array}{l}
\widehat{\mathbb{E}}_{q_k}[\Gamma(x)]
=\sum_{i=1}^{N_k}\widehat{w}_k^i\Gamma(x_k^i),\\
\widehat{Var}_{\theta_k}[\Gamma(x)]=\frac{1}{N_k-1}
\sum_{i=1}^{N_k}\Gamma(x_k^i)\Gamma(x_k^i)^T-
\frac{1}{N^2_k-N_k}\left(\sum_{i=1}^{N_k}\Gamma(x_k^i)\right)
\left(\sum_{i=1}^{N_k}\Gamma(x_k^i)\right)^T.
\end{array}\right.
\end{equation*}
Estimate the gradient $g_k$ by
$\widehat{g}_k :=\widehat{\mathbb{E}}_{q_k}[\Gamma(x)]
-\mathbb{E}_{\theta_k}
\left[\Gamma(x)\right].$
\\
4. \textbf{Updating}: Update the sampling distribution parameter $\theta$ according to
\begin{equation*}
\theta_{k+1}=\Pi_{\widetilde{\Theta}}\left\{\theta_k+\beta_k\left(
\widehat{Var}_{\theta_k}[\Gamma(x)]+\epsilon I\right)^{-1}\widehat{g}_k\right\},
\end{equation*}
where $\widetilde{\Theta}\subseteq \Theta$ is a non-empty compact and convex constraint set.
\\
5. \textbf{Stopping}: Check if some stopping criterion is satisfied. If yes, stop and return the current best sampled solution; else, set $k:= k+1$ and go back to step 2.
\end{algorithm}

In the initialization step (step 1) of GASS-CVaR, the conditions on the sample size and step size sequences are imposed to facilitate the convergence of the algorithm. They are typical requirements for a stochastic approximation algorithm.
Since in the sampling step (step 2) the CVaR values are estimated, the convergence of the original GASS algorithm, which is designed for deterministic optimization, does not directly apply to GASS-CVaR.
We will show the convergence of GASS-CVaR later.
In the estimation step (step 3), as mentioned before, one common choice of the shape function $S_\theta(\cdot)$ is the one in (\ref{eq.2.2}). Moreover, the quantile level $\rho$ in (\ref{eq.2.3}) controls the percentile of elite samples that are used to update the sampling distribution at the next iteration, and balances between the exploitation of the neighborhood of current best solutions and the exploration of the entire solution space. For example, when a smaller $\rho$ is used, less elite samples are used in the updating of the sampling distribution, and thus less emphasis is put on exploration.
In the updating step (step 4), the iterate is projected onto a convex and compact subset $\widetilde{\Theta} \subseteq \Theta$ instead of $\Theta$, in order to guarantee numerical stability and fast computation of the projection.
In the stopping step (step 5), a common stopping criterion used in practice is that the norm of the gradient falls below a pre-specified threshold.

\subsection{GASS with Adaptive Risk Levels} \label{sec3.subsec1:Adaptive}

When the risk level of interest $\alpha^\ast$ is close to $1$, implementing GASS-CVaR could be computationally expensive, since in step 2 the CVaR evaluation requires a large sample size $M_k$ to obtain a good CVaR estimator. This issue is more severe as $\alpha^\ast$ gets closer to $1$.
For example, for a fixed $x$, suppose we want to estimate $C_{\alpha}(x)$ at three different risk levels: $\alpha_1=0$, $\alpha_2=0.90$, and $\alpha_3=0.99$, where we note that $C_{\alpha_1=0}(x)=\mathbb{E}[l(x,\xi_x)]$ is the expected loss.
Loosely speaking, to achieve the same level of accuracy in CVaR estimation, the corresponding sample sizes $M_1$, $M_2$, and $M_3$ should result in equal ``effective'' sample sizes.
In particular, theoretically using $M_i$ i.i.d. samples to estimate $C_{\alpha_i}(x)$ results in $(1-\alpha_i)M_i$ effective samples, since the rest $\alpha_iM_i$ samples result in a value of zero. This implies that
$(1-\alpha_1)M_1=(1-\alpha_2)M_2=(1-\alpha_3)M_3$ for equal effective sample sizes.
Thus, $M_2=(1-\alpha_1)/(1-\alpha_2)\cdot M_1=10 \cdot M_1$ and $M_3=(1-\alpha_1)/(1-\alpha_3)\cdot M_1=100 \cdot M_1$.
Therefore, the sample size required for accurate CVaR estimation could be easily up to tens of times even hundreds of times compared with the sample size required for accurate estimation of expectation.

To save simulation budget and improve the overall efficiency of GASS-CVaR, we propose to initialize the algorithm at a small risk level $\alpha_0$ (e.g., $\alpha_0=0$), and adaptively increase the risk level $\alpha_k$ at every iteration until the target risk level $\alpha^{\ast}$ is achieved while the algorithm converges at the same time.
Since a lower risk level implies that a smaller $M_k$ is required to achieve the desired accuracy for CVaR estimation, the hope is to adaptively save simulation budget at each iteration by solving a problem that is similar to the original one but less computationally expensive.

A good updating rule on the risk level should
1) achieve significant budget savings when the algorithm is in the ``warm-up'' phase, i.e., when it puts more emphasis on the exploration of the entire solution space;
2) solve problems that are close to the original one when the algorithm is in the ``convergence'' phase, i.e., when it puts more emphasis on the exploitation of the promising region that has been identified.
The key to such an updating rule lies in finding an empirical signal on the algorithm's emphasis between exploration and exploitation.

Note that GASS-CVaR maintains the structure of a gradient-based optimization scheme, and thus the gradient $g_k$ (even the Hessian) used in the updating rule of sampling distribution parameter could be regarded as an empirical signal on the algorithm's balance between exploration and exploitation.
Loosely speaking, when the norm of $g_k$ is relatively large, the sampling distribution parameter at next iteration, $\theta_{k+1}$, will differ from $\theta_{k}$ significantly.
This means the algorithm is in the ``warm-up'' phase, where different regions of the solution space are being explored.
When the norm of $g_k$ is small, $\theta_{k+1}$ is expected to be close to $\theta_{k}$.
This means the algorithm is in the ``convergence'' phase, where an identified promising region is being exploited.
Therefore, it is natural to design the updating rule on risk level using the information contained in the gradient $g_k$ obtained at every iteration.
For example, note that GASS-CVaR converges when the norm of $g_k$ hits zero. Then naturally one could increase the risk level at every iteration proportionally to the decrease in the norm of $g_k$ from previous iteration, which ensures that the target risk level $\alpha^\ast$ is achieved when the gradient hits zero, i.e., when the algorithm converges.

In particular, we propose an updating scheme on the risk level as follows.
\begin{equation}\label{eq.3.1}
\alpha_{k+1}=\;\left\{\begin{array}{ll}
\alpha^\ast-\frac{\enorm{g_k}}{\enorm{g_{k-1}}}
\left(\alpha^\ast-\alpha_k\right), &\mbox{if}\; \enorm{g_k}< \enorm{g_{k-1}},\\
\alpha_k, &o/w,
\end{array}\right.
\end{equation}
where $\enorm{\cdot}$ is the vector Euclidean norm.
Note that the updating rule (\ref{eq.3.1}) ensures that $\alpha_k$ is non-decreasing and bounded above by $\alpha^\ast$, with the hope that $\alpha_k$ will eventually converge to $\alpha^\ast$.
Furthermore, when $\enorm{{g}_k}<\enorm{{g}_{k-1}}$, we can rewrite (\ref{eq.3.1}) as $\frac{\alpha^\ast-\alpha_{k+1}}
{\alpha^\ast-\alpha_{k}}=\frac{\enorm{{g}_k}} {\enorm{{g}_{k-1}}}$.
Loosely speaking, it implies the increase in the risk level for next iteration is proportional to the decrease in the norm of the gradient from previous iteration. It also ensures that $\alpha^\ast$ is achieved when the norm of the gradient hits zero, i.e., when the algorithm converges. We do point out that more sophisticated updating rules on the risk level could be incorporated in the future.

Now we present the following Algorithm \ref{alg.3.2}, which is referred to as GASS-CVaR-ARL, for simulation optimization of CVaR with adaptive risk levels.

\begin{algorithm}
\caption{\textbf{GASS-CVaR} with \textbf{A}daptive \textbf{R}isk \textbf{L}evels}\label{alg.3.2}
1. \textbf{Initialization}: Initialize the algorithm similar to step 1 of GASS-CVaR. Set initial risk level $\alpha_0=0$.
\\
2. \textbf{Sampling}:
Draw candidate solutions and simulate the loss distribution scenarios same as step 2 of GASS-CVaR.
Estimate ${C}_{\alpha_k}(x_k^i)$, CVaR of the loss at the risk level $\alpha_k$, by
$$
\widehat{C}_{\alpha_k}(x_k^i)=l\left(x_k^i, \xi_k^{i,(\lceil \alpha_k M_k\rceil)}\right)+\frac{1}{M_k(1-\alpha_k)}
\sum_{j=1}^{M_k}\left(l\left(x_k^i, \xi_k^{i,j}\right)-
l\left(x_k^i, \xi_k^{i,(\lceil \alpha_k M_k\rceil)}\right)\right)^+.
$$
Record the best candidate solution $x_k^\ast$ found at this iteration:
$x_k^\ast=\arg\min_i \;\widehat{C}_{\alpha_k}(x_k^i)$.
\\
3. \textbf{Estimation}: Compute the normalized weights $\widebar{w}_k^i$ according to
$$
\widebar{w}_k^i=\frac{\widehat{S}_{\theta_k}
\left(-\widehat{C}_{\alpha_k}(x_k^i)\right)}
{\sum_{j=1}^{N_k}\widehat{S}_{\theta_k}
\left(-\widehat{C}_{\alpha_k}(x_k^j)\right)}, \quad i=1,...,N_k,
$$
and then estimate $\mathbb{E}_{q_k}[\Gamma(x)]$ and $Var_{\theta_k}[\Gamma(x)]$  via
\begin{equation*}
\left\{\begin{array}{l}
\widebar{\mathbb{E}}_{q_k}[\Gamma(x)]
=\sum_{i=1}^{N_k}\widebar{w}_k^i\Gamma(x_k^i),\\
\widehat{Var}_{\theta_k}[\Gamma(x)]=\frac{1}{N_k-1}
\sum_{i=1}^{N_k}\Gamma(x_k^i)\Gamma(x_k^i)^T-
\frac{1}{N^2_k-N_k}\left(\sum_{i=1}^{N_k}\Gamma(x_k^i)\right)
\left(\sum_{i=1}^{N_k}\Gamma(x_k^i)\right)^T.
\end{array}\right.
\end{equation*}
Estimate the gradient ${g}_k$ by
$\widebar{g}_k :=\widebar{\mathbb{E}}_{q_k}[\Gamma(x)]
-\mathbb{E}_{\theta_k}
\left[\Gamma(x)\right].$
\\
4. \textbf{Updating}: Update the sampling distribution parameter $\theta$ according to
\begin{equation*}
\theta_{k+1}=\Pi_{\widetilde{\Theta}}\left\{\theta_k+\beta_k\left(
\widehat{Var}_{\theta_k}[\Gamma(x)]+\epsilon I\right)^{-1}\widebar{g}_k\right\},
\end{equation*}
where $\widetilde{\Theta}\subseteq \Theta$ is a non-empty compact and convex constraint set; then update the risk level $\alpha$ according to
\begin{equation}\label{eq.3.2}
\alpha_{k+1}=\;\left\{\begin{array}{ll}
\alpha^\ast-\frac{\enorm{\widebar{g}_k}}{\enorm{\widebar{g}_{k-1}}}
\left(\alpha^\ast-\alpha_k\right), &\mbox{if}\; \enorm{\widebar{g}_k}< \enorm{\widebar{g}_{k-1}},\\
\alpha_k, &o/w.
\end{array}\right.
\end{equation}
\\
5. \textbf{Stopping}: Check if some stopping criterion is satisfied. If yes, stop and return
$x^{\ast}=\arg\min_k \;\widehat{C}_{\alpha^\ast}(x_k^\ast)$ and $\widehat{C}_{\alpha^\ast}(x^\ast)$ via simulation; else, set $k:= k+1$ and go back to step 2.
\end{algorithm}

In the sampling step (step 2) of GASS-CVaR-ARL, since the current risk level $\alpha_k$ is smaller than the target risk level $\alpha^\ast$, we could use a sample size $M_k$ smaller than the one used in GASS-CVaR to estimate the CVaR values at risk level $\alpha_k$. For example, suppose one wants to keep the ``effective'' sample size $(1-\alpha_k)M_k$ as a constant. Then, in the initial iterations of the algorithm the budget savings can be up to tens of times even hundreds of times (equal to $(1-\alpha_k)/(1-\alpha^\ast)$ precisely) since $\alpha_k$ is close to $\alpha_0=0$.
The best candidate solution generated at each iteration is also recorded, where note that at $k^{th}$ iteration it is identified by the minimum CVaR value at risk level $\alpha_k$. So it is a good solution to the CVaR minimization problem as if the target risk level is $\alpha_k$.

In the estimation step (step 3), note that the estimation of $\{w\}$, $\mathbb{E}_{q_k}[\Gamma(x)]$, and thus $g_k$ differs from the estimation in GASS-CVaR, since $\alpha_k$ instead of $\alpha^\ast$ is used as the risk level at the $k^{th}$ iteration.
One could also view $\widebar{g}_k$ as an approximation of the gradient for the reformulated problem of $\max_{x\in\mathcal{X}} -C_{\alpha_k}(x)$.
This implies that at each step GASS-CVaR-ARL solves a CVaR optimization problem that is structurally similar to the original one but less computationally intensive.

In the updating step (step 4), the updating rule (\ref{eq.3.2})
is an implementable version of (\ref{eq.3.1}), with the gradient $g_k$ being replaced by $\widebar{g}_k$.
Note that it still ensures that $\alpha_k$ is non-decreasing bounded above by $\alpha^\ast$, and $\alpha^\ast$ is achieved when the norm of $\widebar{g}_k$ hits zero, i.e., when the algorithm converges.

In the stopping step (step 5), finding the best solution to the original CVaR optimization problem is achieved via evaluating and comparing the CVaR values at the target risk level $\alpha^\ast$ for all the best candidate solutions found so far, and thus additional simulation budget is required; however, it is insignificant compared with the overall budget consumed.

Recall that, in GASS-CVaR-ARL, the risk level used at each iteration is updated in accordance with the decrease in the norm of the gradient.
It implies that the updating rule (\ref{eq.3.1}) keeps track of the algorithm's balance between the exploration of the entire solution space and the exploitation of an identified promising region, and then makes adjustments on the risk level accordingly.
Therefore, in the ``warm-up'' phase of the algorithm, using a small risk level has little negative effect on the algorithm progress since the algorithm puts most of its emphasis on exploration; in the ``convergence'' phase of the algorithm, the risk level $\alpha_k$ is close to $\alpha^\ast$, and essentially the algorithm is solving problems that are very close to the original one.
Thus, intuitively, we expect the number of iterations that GASS-CVaR-ARL takes to converge to be similar to the one that GASS-CVaR takes to converge, which is also verified by the numerical tests presented in Section \ref{sec5:Numerical}.
Since GASS-CVaR-ARL saves simulation budget at every iteration, total budget saving is achieved.

\section{Convergence Analysis}\label{sec4:Convergence}

Let us first analyze the convergence properties of GASS-CVaR (Algorithm \ref{alg.3.1}). The analysis will rely mainly on the convergence analysis of GASS in \cite{zhou2014gradient} as well as the classic results in stochastic approximation methods and algorithms (see, e.g., \cite{kushner2003stochastic}, \cite{borkarstochastic}, \cite{kushner2010stochastic}, and \cite{kushner2012stochastic}).
The main idea is to reformulate the updating scheme on $\theta_k$ in GASS-CVaR as a generalized Robbins-Monro recursive algorithm in solving a constrained ordinary differential equation (ODE) of $\theta$, and then show the corresponding bias term and noise term in the reformulated updating scheme are bounded in appropriate asymptotical sense so that the sequence $\{\theta_k\}$ generated by the updating scheme converges to a limit set of the ODE w.p.1.

Following the above road map, let us first reformulate the parameter updating scheme in GASS-CVaR
\begin{equation}\label{eq.4.1}
\begin{split}
\theta_{k+1}&=\Pi_{\widetilde{\Theta}}\left\{\theta_k+\beta_k\left(
\widehat{Var}_{\theta_k}[\Gamma(x)]+\epsilon I\right)^{-1}\widehat{g}_k\right\}\\
&=\Pi_{\widetilde{\Theta}}\left\{\theta_k+\beta_k\left(
\widehat{Var}_{\theta_k}[\Gamma(x)]+\epsilon I\right)^{-1}\left(
\widehat{\mathbb{E}}_{q_k}[\Gamma(x)]
-\mathbb{E}_{\theta_k}\left[\Gamma(x)\right]\right)\right\},
\end{split}
\end{equation}
as
\begin{equation}\label{eq.4.2}
\theta_{k+1}=\theta_k+\beta_k\left[G(\theta_k)+b_k+e_k+p_k\right].
\end{equation}
Here
\begin{eqnarray*}
G(\theta_k)&\overset{\triangle}=&V_k^{-1}\left(
\mathbb{E}_{q_k}[\Gamma(x)]
-\mathbb{E}_{\theta_k}\left[\Gamma(x)\right]\right),\\
b_k&\overset{\triangle}=&\widehat{V}_k^{-1}\left(
\widehat{\mathbb{E}}_{q_k}[\Gamma(x)]-
\widetilde{\mathbb{E}}_{q_k}[\Gamma(x)]\right),\\
e_k&\overset{\triangle}=&\left(\widehat{V}_k^{-1}-V_k^{-1}\right)
\left(\widetilde{\mathbb{E}}_{q_k}[\Gamma(x)]
-\mathbb{E}_{\theta_k}\left[\Gamma(x)\right]\right)+
V_k^{-1}\left(
\widetilde{\mathbb{E}}_{q_k}[\Gamma(x)]
-\mathbb{E}_{q_k}\left[\Gamma(x)\right]\right),
\end{eqnarray*}
and $p_k$ is the resulted projection error term, where for simplicity we denote
\begin{equation*}
V_k\overset{\triangle}=\left(
Var_{\theta_k}[\Gamma(x)]+\epsilon I\right)\quad\mbox{and}\quad
\widehat{V}_k\overset{\triangle}=\left(
\widehat{Var}_{\theta_k}[\Gamma(x)]+\epsilon I\right).
\end{equation*}
In (\ref{eq.4.2}) the term $G(\theta_k)$ is the gradient vector field in a standard stochastic approximation algorithm, the term $b_k$ represents the bias in estimating $\widetilde{\mathbb{E}}_{q_k}[\Gamma(x)]$ caused by the inexact evaluation of the shape function, the term $e_k$ represents the simulation noise in the estimators $\widehat{Var}_{\theta_k}[\Gamma(x)]$ and $\widetilde{\mathbb{E}}_{q_k}[\Gamma(x)]$, and the term $p_k$ represents the projection error after taking the current iterate back onto the constraint set $\widetilde{\Theta}$ with minimum Euclidean norm.
Note that the bias term $b_k$ is caused by both the outer-layer sampling on the solution space and the inner-layer simulation of the loss distribution;
however, the noise term $e_k$ accounts for the error due to the outer-layer sampling only, since both $\widehat{Var}_{\theta_k}[\Gamma(x)]$ and
$$
\widetilde{\mathbb{E}}_{q_k}[\Gamma(x)]=
\sum_{i=1}^{N_k}w_k^i \Gamma(x^i_k),
\quad \mbox{where}\quad
w_k^i=\frac{S_{\theta_k}(-C_{\alpha^\ast}(x_k^i))}
{\sum_{j=1}^{N_k}S_{\theta_k}(-C_{\alpha^\ast}(x_k^i))}
$$
do not involve the inner-layer sampling of $l(x;\xi_x)$.

Now let us introduce the assumptions on the algorithm and $l(x;\xi_x)$ for the convergence of the algorithm.
The following set of assumptions is on the algorithm parameters and the choice of the exponential family of densities.
It largely follows from the standard assumptions for a generalized stochastic approximation algorithm.
\begin{assumption}\label{asm.4.1}~
\begin{itemize}
\item[(i)] The step size sequence $\{\beta_k\}$ satisfies that $\beta_k>0$ for all $k$, $\beta_k\searrow0$ as $k\rightarrow\infty$, $\sum_{k=0}^{\infty}\beta_k=\infty$ and $\sum_{k=0}^{\infty}\beta^2_k<\infty$.

\item[(ii)] The outer-layer sample size sequence $\{N_k\}$ satisfies $N_k=N_0\cdot k^\tau$ for some constant $\tau> 0$. Furthermore, the sequences $\{\beta_k\}$ and $\{N_k\}$ jointly satisfies $\frac{\beta_k}{\sqrt{N_k}}=O(k^{-\zeta})$ for some constant $\zeta>1$.

\item[(iii)] The inner-layer sample size sequence $M_k$ satisfies that $M_k\nearrow \infty$ as $k\rightarrow \infty$.

\item[(iv)] The sufficient statistics $\Gamma(x)$ of the chosen exponential family of densities is bounded on $\mathcal{X}$.
\end{itemize}
\end{assumption}

In the above set of assumptions, Assumption \ref{asm.4.1}.(i) follows from the typical step size assumption in a gradient-based optimization algorithm.
Assumption \ref{asm.4.1}.(ii) ensures that the outer-layer sample size $N_k$ increases to infinity no slower than certain speed given a choice of the step size sequence, and it can be easily satisfied.
For example, if $\beta_k=O(1/k)$, then $N_k=N_0\cdot k^\tau$ for an arbitrary constant $\tau> 0$ is sufficient for $\frac{\beta_k}{\sqrt{N_k}}=O(k^{-\zeta})$ to hold for some constant $\zeta>1$.
Assumption \ref{asm.4.1}.(iii) ensures that the error of the CVaR estimators caused by the inner-layer simulation of the loss distribution vanishes as $k\rightarrow \infty$.
Assumption \ref{asm.4.1}.(iv) is to bound the expectation and variance terms of the sufficient statistics in the algorithm. It holds for many exponential families used in practice.
For example, when the solution space $\mathcal{X}$ is a nonempty compact set, the continuity of the function $\Gamma(\cdot)$ will be sufficient for Assumption \ref{asm.4.1}.(iv) to hold.

The next set of assumptions is on the regularity conditions of the loss function $l(x,\xi_x)$.
As noted previously, the bias term $b_k$ is caused by the inexact evaluation of the shape function $S_{\theta_k}(\cdot)$.
When $S_{\theta_k}(\cdot)$ takes the form of (\ref{eq.2.2}), $b_k$ is caused by the error in estimating the $(1-\rho)$-quantile $\gamma_{\theta_k}$ in (\ref{eq.2.3}) as well as the error in Monte Carlo estimation of the CVaR values.
Specifically, recall that for a fixed $x$,
\begin{equation}\label{eq.4.3}
\widehat{S}_{\theta_k}\left(-\widehat{C}_{\alpha^\ast}(x)\right)=
\frac{1}{1+\exp\left(-S_o\left(-\widehat{C}_{\alpha^\ast}(x)
-\widehat{\gamma}_{\theta_k}
\right)\right)},
\end{equation}
where $\widehat{C}_{\alpha^\ast}(x)$ is the CVaR estimator given in (\ref{eq.2.10}), and $\widehat{\gamma}_{\theta_k}$ is the sample $(1-\rho)$-quantile of $\{-\widehat{C}_{\alpha^\ast}(x_k^i): i=1,...,N_k\}$, i.e., $\widehat{\gamma}_{\theta_k}$ is the $(\lceil(1-\rho)N_k\rceil)^{th}$ order statistic of $\{-\widehat{C}_{\alpha^\ast}(x_k^i): i=1,...,N_k\}$.
Since $\gamma_{\theta_k}$ could be viewed as the $(1-\rho)$-level Value-at-Risk (VaR) of $-C_{\alpha^\ast}(x)$ w.r.t. the sampling distribution $f(x;\theta_k)$, then $\widehat{\gamma}_{\theta_k}$ could be regarded as a nested risk estimator in which the outer-layer simulation is on estimation of VaR and the inner-layer is on estimation of CVaR.
Hence, bounding the bias term $b_k$ reduces to bounding the errors of the nested risk estimator $\widehat{\gamma}_{\theta_k}$ as well as the one-layer CVaR estimator $\widehat{C}_{\alpha^\ast}(x)$. Here we will resort to the asymptotic analysis of nested risk estimators in \cite{gordy2010nested}, and \cite{zhu2016risk}.

To this end, let us rewrite the CVaR estimator $\widehat{C}_{\alpha^\ast}(x)$ in (\ref{eq.2.10}) as
\begin{equation*}
\widehat{C}_{\alpha^\ast}(x)
=C_{\alpha^\ast}(x)+\frac{1}{\sqrt{M_k}}\cdot\mathcal{E}_k(x),\quad \forall x\in\mathcal{X},
\end{equation*}
where $\mathcal{E}_k(x)$ is the standardized error of the CVaR estimator.
Note that by the asymptotic normality of $\widehat{C}_{\alpha^\ast}(x)$, under appropriate regularity conditions $\mathcal{E}_k(x)$ has a limiting distribution as $k\rightarrow\infty$.
Thus, the effect of the diminishing noise term $\mathcal{E}_k(x)/\sqrt{M_k}$ on the distribution of $\widehat{C}_{\alpha^\ast}(x)$ will vanish as $M\rightarrow \infty$. Hence, we expect the ``distance'' between the distribution of $\widehat{C}_{\alpha^\ast}(x)$ and the distribution of $C_{\alpha^\ast}(x)$ to vanish as $M\rightarrow \infty$. That is,
the p.d.f. of $\widehat{C}_{\alpha^\ast}(x)$ converges to the p.d.f. of $C_{\alpha^\ast}(x)$.
The following set of assumptions, which is referred to Assumption \ref{asm.4.2}, guarantees that the convergence of the p.d.f. is sufficiently fast. It largely follows from Assumption 1 in \cite{gordy2010nested} and Assumption 3.2 in \cite{zhu2016risk}.

\begin{assumption}\label{asm.4.2}~
\begin{itemize}

\item[(i)] For all $x\in \mathcal{X}$, the loss distribution $l(x,\xi_x)$ has finite second moment; moreover, for all $\theta\in {\Theta}$, the CVaR function $C_{\alpha^\ast}(x)$, which is a random variable under the distribution $f(\cdot;\theta)$, has finite second moment.

\item[(ii)] For all $\theta\in {\Theta}$ and each $k$, the joint density $d_k(c, e)$ of ($C_{\alpha^\ast}(x), \mathcal{E}_k(x))$, and its partial derivatives $\frac{\partial}{\partial c}d_k(c, e)$ and $\frac{\partial^2}{\partial c^2}d_k(c, e)$ exist for all pairs of $(c, e)$.

\item[(iii)] For all $\theta\in {\Theta}$ and each $k$, there exist nonnegative functions $D_{0,k}(\cdot)$, $D_{1,k}(\cdot)$ and $D_{2,k}(\cdot)$ such that $d_k(c, e)<D_{0,k}(e)$, $\frac{\partial}{\partial c}d_k(c, e)<D_{1,k}(e)$, and $\frac{\partial^2}{\partial c^2}d_k(c, e)<D_{2,k}(e)$ for all $(c, e)$. Furthermore, for all $\theta\in {\Theta}$, $\sup_{k}\int|e|^r D_{i,k}(e)de<\infty$ for $i=0,1,2$, and $0\le r\le 4$.

\end{itemize}
\end{assumption}
In the above assumption, Assumption \ref{asm.4.2}.(i) ensures that a one-layer VaR or CVaR estimator defined in (\ref{eq.2.10}) or (\ref{eq.2.11}) is strongly consistent and asymptotically normally distributed, and thus the standardized estimation error $\mathcal{E}_k(x)$ has a limiting distribution as $k\rightarrow\infty$.
Assumption \ref{asm.4.2}.(ii) and \ref{asm.4.2}.(iii) further ensure that the p.d.f. of $\widehat{C}_{\alpha^\ast}(x)$ converges to the p.d.f. of $C_{\alpha^\ast}(x)$ sufficiently fast.
This will imply the strong consistency of the nested risk estimator $\widehat{\gamma}_{\theta_k}$ and further the convergence of the approximate shape function $\widehat{S}_{\theta_k}\left(-\widehat{C}_{\alpha^\ast}(x)\right)$, as presented in the following Lemma \ref{lem.4.1}.


\begin{lemma}\label{lem.4.1}
Suppose the shape function $S_{\theta_k}(\cdot)$ takes the form
\begin{equation*}
S_{\theta_k}\left(-C_{\alpha^\ast}(x)\right)=
\frac{1}{1+\exp\left(-S_o\left(-C_{\alpha^\ast}(x)-\gamma_{\theta_k}
\right)\right)},
\end{equation*}
where $S_o$ is a large positive constant, and
$\gamma_{\theta_k}\overset{\triangle}=\sup_{r}\left\{
r: P_{f(\cdot;\theta_k)}\left\{x\in \mathcal{X}: -C_{\alpha^\ast}(x)\ge r\right\} \ge \rho\right\}$ is the $(1-\rho)$-quantile of $(-C_{\alpha^\ast}(x))$ w.r.t. $f(\cdot;\theta_k)$. Further suppose that  $S_{\theta_k}\left(-C_{\alpha^\ast}(x)\right)$ is approximated by $\widehat{S}_{\theta_k}\left(-\widehat{C}_{\alpha^\ast}(x)\right)$ as in (\ref{eq.4.3}). Then under Assumption \ref{asm.4.1}.(ii), \ref{asm.4.1}.(iii) and Assumption \ref{asm.4.2}, we have
\begin{equation}\label{eq.4.4}
\lim\limits_{k\rightarrow\infty}\left|
\widehat{S}_{\theta_k}\left(-\widehat{C}_{\alpha^\ast}(x)\right)
-S_{\theta_k}\left(-C_{\alpha^\ast}(x)\right)\right|=0,\quad
w.p.1,\quad \forall x\in\mathcal{X}.
\end{equation}
\end{lemma}

The main idea of the proof is to show $\widehat{C}_{\alpha^\ast}(x)\rightarrow {C}_{\alpha^\ast}(x)$ w.p.1 and $\widehat{\gamma}_{\theta_k}\rightarrow\gamma_{\theta_k}$ w.p.1 as $k\rightarrow \infty$.
The detailed proof is included in the appendix.
Following the road map and based on Lemma 1, we next show that the bias term $b_k$ converges to zero w.p.1. as $k\rightarrow\infty$, as presented in Lemma \ref{lem.4.2} below.

\begin{lemma}\label{lem.4.2}
Under Assumption \ref{asm.4.1} and Assumption \ref{asm.4.2}, we have \begin{equation}\label{eq.4.5}
\lim_{k\rightarrow \infty} \enorm{b_k}=0, \; w.p.1,
\end{equation}
where recall that $\enorm{b_k}$ is the vector Euclidean norm of $b_k$.
\end{lemma}

The proof of Lemma \ref{lem.4.2} is included in the appendix. Continuing the road map, we next show that the summed tail error goes to zero w.p.1, as presented in the following Lemma \ref{lem.4.3}.

\begin{lemma}\label{lem.4.3}
Under Assumption \ref{asm.4.1}, we have
\begin{equation}\label{eq.4.6}
\lim_{k\rightarrow\infty} \left\{
\sup_{n: 0\le \sum_{i=k}^{n-1}\beta_i\le T}\enorm{\sum_{i=k}^n \beta_i e_i}\right\}=0,\quad w.p.1
\end{equation}
for all $T\ge0$.
\end{lemma}

Lemma \ref{lem.4.3} is identical to Lemma 2 in \cite{zhou2014gradient}, so we omit the proof here. With the above lemmas, we now proceed to the main result on the convergence of Algorithm \ref{alg.3.1}.

Given an arbitrary $\theta\in \widetilde{\Theta}$, a set $\mathcal{C}(\theta)$ is defined as follows. For $\theta$ that lies in the interior of $\widetilde{\Theta}$, let $\mathcal{C}(\theta)=\{0\}$; for $\theta$ that lies on the boundary of $\widetilde{\Theta}$, let $\mathcal{C}(\theta)$ be the infinite convex cone generated by the outer normals at $\theta$ of the faces on which $\theta$ lies (see, e.g., \cite{kushner2010stochastic} pp. 89). Then the updating scheme (\ref{eq.4.2}) in GASS-CVaR could be viewed as a noisy discretization of a constrained ODE for $\{\theta(t): t\ge0\}$:
\begin{equation}\label{eq.4.7}
\dot{\theta}(t)=G(\theta(t))+p(t), \;p(t)\in -\mathcal{C}(\theta(t)), \; t\ge 0,
\end{equation}
where $p(t)$ is the minimum force to take $\theta(t)$ back to the set $\widetilde{\Theta}$.
Using the ODE approach for the convergence of the Robbins-Monro Algorithm (see, e.g., \cite{kushner2010stochastic}), we can show that the sequence $\{\theta_k\}$ generated by (\ref{eq.4.1}) converges to a limit set of the ODE (\ref{eq.4.7}). In particular, we have the following theorem.

\begin{theorem}\label{thm.4.1}{\textbf{Convergence of GASS-CVaR}.}
Suppose Assumption \ref{asm.4.1} and Assumption \ref{asm.4.2} hold. Then the sequence $\{\theta_k\}$ generated by (\ref{eq.4.1}) converges to a limit set of the ODE (\ref{eq.4.7}) w.p.1. Furthermore, if the limit sets of (\ref{eq.4.7}) are isolated equilibrium points, then $\{\theta_k\}$ converges to a unique equilibrium point w.p.1.
\end{theorem}

Theorem \ref{thm.4.1} is a direct consequence of Theorem 2 in \cite{kushner2010stochastic} with Lemma \ref{lem.4.2} and Lemma \ref{lem.4.3} above.
Starting with the convergence of GASS-CVaR, we will show the convergence of the algorithm GASS-CVaR-ARL. The intuition is as follows.

Recall that the updating scheme on $\theta$ in GASS-CVaR-ARL is
\begin{equation}\label{eq.4.8}
\begin{split}
\theta_{k+1}&=\Pi_{\widetilde{\Theta}}\left\{\theta_k+\beta_k\left(
\widehat{Var}_{\theta_k}[\Gamma(x)]+\epsilon I\right)^{-1}\widebar{g}_k\right\}\\
&=\Pi_{\widetilde{\Theta}}\left\{\theta_k+\beta_k\left(
\widehat{Var}_{\theta_k}[\Gamma(x)]+\epsilon I\right)^{-1}\left(\widebar{\mathbb{E}}_{q_k}[\Gamma(x)]
-\mathbb{E}_{\theta_k}
\left[\Gamma(x)\right]\right)\right\}.
\end{split}
\end{equation}
Compared with the updating scheme (\ref{eq.4.1}) on $\theta$ in  GASS-CVaR, we could see that the approximate expectation term $\widehat{\mathbb{E}}_{q_k}[\Gamma(x)]$ in (\ref{eq.4.1}) is replaced by $\widebar{\mathbb{E}}_{q_k}[\Gamma(x)]$ in (\ref{eq.4.8}) in estimating the gradient $g_k$. Note that the updating scheme for the risk level $\alpha_k$ in (\ref{eq.3.2}) guarantees that $\alpha_k$ is non-decreasing and bounded above by the target risk level $\alpha^{\ast}$. Thus, the limit of the risk level sequence $\{\alpha_k\}$ exists. If we are able to show that the limit is $\alpha^\ast$, then the difference between $\widebar{\mathbb{E}}_{q_k}[\Gamma(x)]$ and $\widehat{\mathbb{E}}_{q_k}[\Gamma(x)]$, i.e., the difference between $\widebar{g}_k$ and $\widehat{g}_k$, will vanish as $k\rightarrow\infty$. The reason is that the normalized weights $\{\widebar{w}_k\}$ in computing $\widebar{\mathbb{E}}_{q_k}[\Gamma(x)]$ will asymptotically approach $\{\widehat{w}_k\}$ in computing $\widehat{\mathbb{E}}_{q_k}[\Gamma(x)]$ as $k\rightarrow\infty$.

Assume by contradiction that $\lim_{k\rightarrow\infty} \alpha_k=\widebar{\alpha}^\ast<\alpha^\ast$. On the one hand, following from above argument, GASS-CVaR-ARL asymptotically approaches GASS-CVaR for the simulation optimization of $C_{\widebar{\alpha}^\ast}(x)$ instead of $C_{\alpha^\ast}(x)$. Therefore, it is convergent, and thus the gradient sequence $\{g_k\}$ approaches zero w.p.1. One the other hand, the sequence $\{\enorm{\widebar{g}_k}\}$ generated by (\ref{eq.3.2}) will always be above a certain positive value w.p.1 (otherwise $\alpha_k$ will converge to $\alpha^\ast$). This contradicts with the fact that $\{g_k\}$ approaches zero w.p.1. We formalize the above analysis in the following Theorem \ref{thm.4.2}. The detailed proof is included in the appendix.
\begin{theorem}\label{thm.4.2}\textbf{Convergence of GASS-CVaR-ARL}.
Suppose Assumption \ref{asm.4.1} and Assumption \ref{asm.4.2} hold. Then the risk level sequence $\{\alpha_k\}$ generated by (\ref{eq.3.2}) converges to the target risk level $\alpha^\ast$ w.p.1, and the sequence $\{\theta_k\}$ generated by (\ref{eq.4.8}) converges to a limit set of the ODE (\ref{eq.4.7}) w.p.1. Furthermore, if the limit sets of (\ref{eq.4.7}) are isolated equilibrium points, then $\{\theta_k\}$ converges to a unique equilibrium point w.p.1.
\end{theorem}

\section{Numerical Experiments}\label{sec5:Numerical}
We carry out numerical tests to compare the performances of GASS-CVaR and GASS-CVaR-ARL.
In particular, the loss functions tested are listed in the following, among which some are designed by adding Gaussian noises to the continuous benchmark functions in \cite{hu2007model}. However, we point out our algorithms do not have much assumption on the structure of the loss function or the noise.
For convenience, let $\mathcal{N}(0,1)$ be a standard one-dimensional Gaussian distribution, and the loss function is in the form of
\begin{equation}\label{eq.5.1}
l_i(x,\xi_x)=L_i(x)+
\left\{\begin{array}{l}
\sqrt{1+100 \sum_{d=1}^D(x_d-1)^2}\cdot \mathcal{N}(0,1),\; i=0, 1, 3, 4,\\
\sqrt{1+100 \sum_{d=1}^D(x_d-2)^2}\cdot \mathcal{N}(0,1),\; i=2, 5,
\end{array}\right.
\end{equation}
where $D$ is the dimension of the solution space.
Specifically, $L_0=\sum_{d=1}^D x_d^2$; $L_1$ and $L_2$ are respectively Powell function and Rosenbrock function, which are badly scaled;
$L_3$ is Rastrigin function, which is multimodal with a large number of local optima;
$L_4$ and $L_5$ are respectively Pint\'{e}r function and Levy function, which are badly-scaled as well as multimodal.
The explicit expressions of $L_i$'s are listed as follows, and we test all functions with $D=10$.
\begin{itemize}
\item[(0)] $L_0(x)=\sum_{d=1}^{D}x_d^2$.

\item[(1)] Powell function $L_1(x)$.
\begin{eqnarray*}
L_1(x)&=&\sum_{d=2}^{D-2}\bigl[(x_{d-1}+10x_d)^2+5(x_{d+1}-x_{d+2})^2+
(x_d-2x_{d+1})^4\\
&&+10(x_{d-1}-x_{d+2})^4\bigr].
\end{eqnarray*}

\item[(2)] Rosenbrock function $L_2(x)$.
\begin{equation*}
L_2(x)=\sum_{d=1}^{D-1}\left[(x_d-1)^2+100(x_d^2-x_{d+1})^2\right].
\end{equation*}

\item[(3)] Rastrigin function $L_3(x)$.
\begin{equation*}
L_3(x)=\sum_{d=1}^{D}(x_d^2-10\cos(2\pi x_d))-10D-1.
\end{equation*}

\item[(4)] Pint\'{e}r function $L_4(x)$.
\begin{eqnarray*}
L_4(x)&=&\Bigl[\sum_{d=1}^{D}dx_d^2+\sum_{d=1}^{D}20d\sin^2(x_{d-1}\sin x_d-x_d+\sin x_{d+1})\\
&&+\sum_{d=1}^{D} d\log_{10}(1+d(x_{d-1}^2-2x_d+3x_{d+1}-\cos x_d+1)^2)\Bigr].
\end{eqnarray*}

\item[(5)] Levy function $L_5(x)$.
\begin{eqnarray*}
L_5(x)&=&-\sin^2(\pi y_1)-\sum_{d=1}^{D-1}\left[
(y_d-1)^2(1+10\sin^2(\pi y_d+1))\right]\\
&&-(y_D-1)^2(1+10\sin^2(2\pi y_D)),
\end{eqnarray*}
where $y_d=1+(x_d-1)/4,~ d=1,...,D$.
\end{itemize}

Note that we add the noise in the above form to make sure the optimal solution becomes different when $\alpha^\ast$ varies. If $C_{\alpha^\ast=0}(x)=L_i(x)=\mathbb{E}_{\xi_x}
[l_i(x, \xi_x)]$ is of interest, then evidently $x^{\circ}=[0,...,0]_d$ is the minimizer for $i=0, 1, 3, 4$, and $x^{\circ}=[1,...,1]_d$ is the minimizer for $i=2, 5$.
As the risk level of interest $\alpha^\ast$ increases, the minimizer of $C_{\alpha^\ast}(x)$, might be very different from $x^{\circ}$. Specifically, the loss distribution of $l_i$ has a relatively large variance at $x^{\circ}$ (note that it has the smallest variance at $x=[x^o_1+1,...,x^o_D+1]$).
This indicates that, as the risk level of interest $\alpha^\ast$ increases, the minimizer of $C_{\alpha^\ast}(x)$ may start to deviate away from $x^{\circ}$ and move towards $x=[x^o_1+1,...,x^o_D+1]$ (this is also verified by our numerical tests), where the loss function is exposed to the lowest amount of noise.
Note that when $\alpha^\ast>0$, except for $l_0$, the minimizers of $C_{\alpha^\ast}(x)$ and the minimum CVaR function values are not analytically available.

In all the implementations, we use independent multivariate normal distribution $\mathcal{N}(\mu_k, \Sigma_k)$ as the parameterized sampling distribution $f(x;\theta_k)$ at iteration $k$, where $\mu_k=(\mu_k^1,...,\mu^D_k)^T$ is the mean parameter and $\Sigma_k=diag((\sigma^1_k)^2,...,(\sigma^D_k)^2)$ is the covariance matrix. Thus,
$
\theta_k=(\mu_k^1,...,\mu^D_k;(\sigma^1_k)^2,...,(\sigma^D_k)^2)^T.
$
The initial mean parameter $\mu_0$ are drawn randomly from the uniform distribution $U[-30,30]^D$, and the initial covariance matrix $\Sigma_0$ is set to be $\Sigma_0=1000I_{D\times D}$, where $I_{D\times D}$ is the identity matrix of dimension $D$.
From the experiment results, we notice that the performance of the algorithms is insensitive to the initial mean parameter as long as the initial covariance matrix is sufficiently large.

At iteration $k$, we use the shape function $S_{\theta_k}(\cdot)$ in the form of expression (\ref{eq.2.2}) with $S_o=10^5$ and $\rho=0.1$ in (\ref{eq.2.3}).
The $(1-\rho)$-quantile $\gamma_{\theta_k}$ is estimated by the $(1-\rho)$ sample quantile of the CVaR estimates for all the candidate solutions generated at this iteration.
The risk level of interest is $\alpha^\ast=0.99$, and in GASS-CVaR-ARL the initial risk level is set to be $\alpha_0=0$.
The sample size of candidate solutions drawn from the sampling distribution is set to be $N_k=1000$, and the sample size used to estimate the CVaR of the loss distribution is set in a way such that the effective sample size is $(1-\alpha_k)M_k=50$.
Therefore, in GASS-CVaR $M_k=50/(1-\alpha^\ast)=50/0.01=5\times 10^3$ for all $k$, and in GASS-CVaR-ARL $M_k=50/(1-\alpha_k)$ at iteration $k$ with initial sample size $M_0=50/(1-0)=50$.
The small positive constant $\epsilon$ used to ensure the positive definiteness of the Hessian is set to be $\epsilon=10^{-10}$, and the step size $\beta_k$ is set to be $\beta_k=50/(k+2000)^{0.6}$, which satisfies the assumptions in step 1 of both two algorithms.

\begin{figure}[htb!]
{
\centering
\begin{tabular}{cc}
\includegraphics[width=0.50\textwidth]{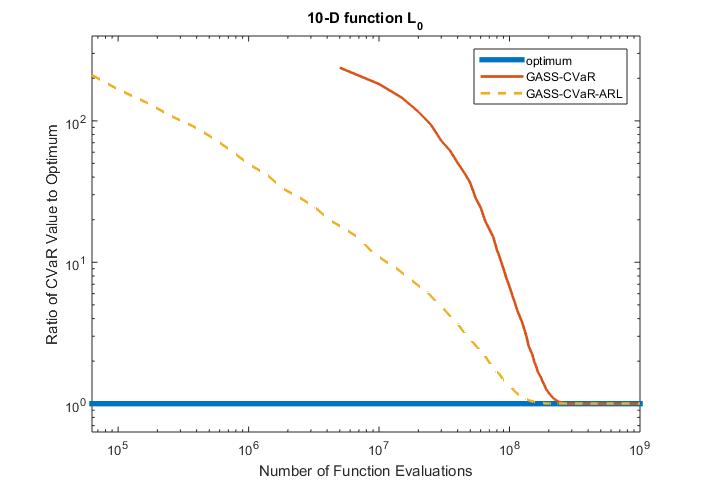}&
\hspace{-5mm}
\includegraphics[width=0.50\textwidth]{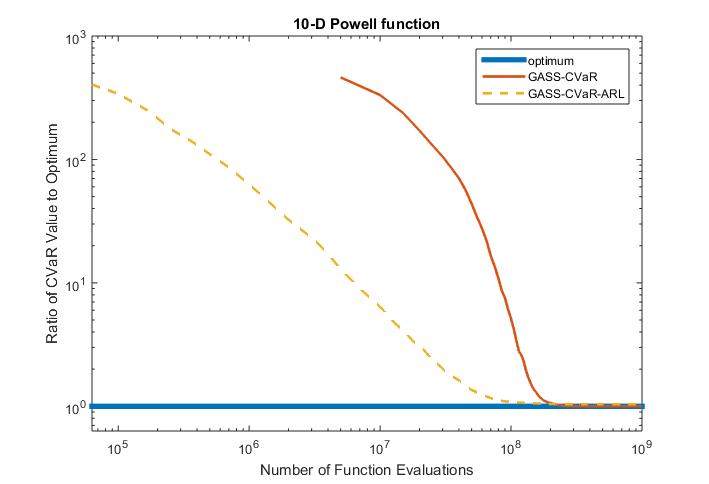}
\\
\includegraphics[width=0.50\textwidth]{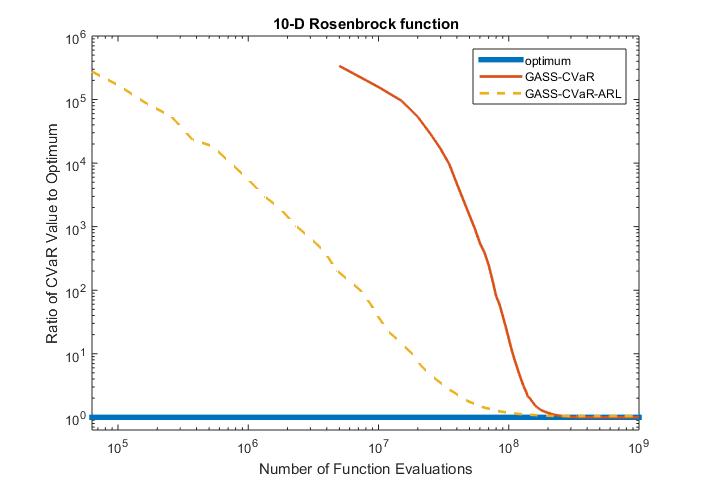}&
\hspace{-5mm}
\includegraphics[width=0.50\textwidth]{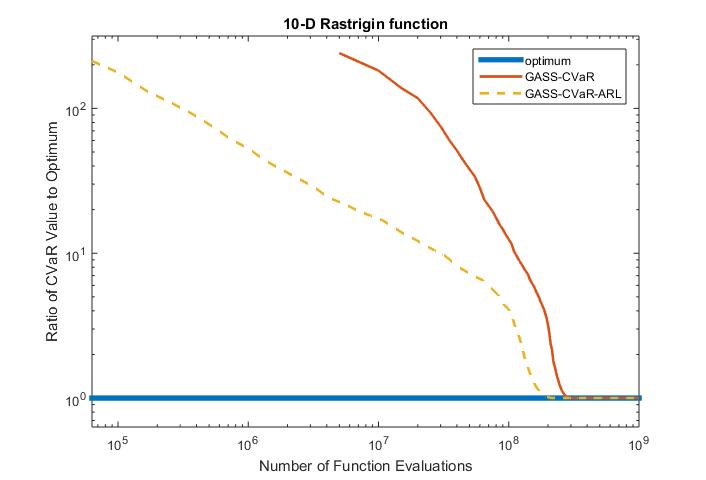}
\\
\includegraphics[width=0.50\textwidth]{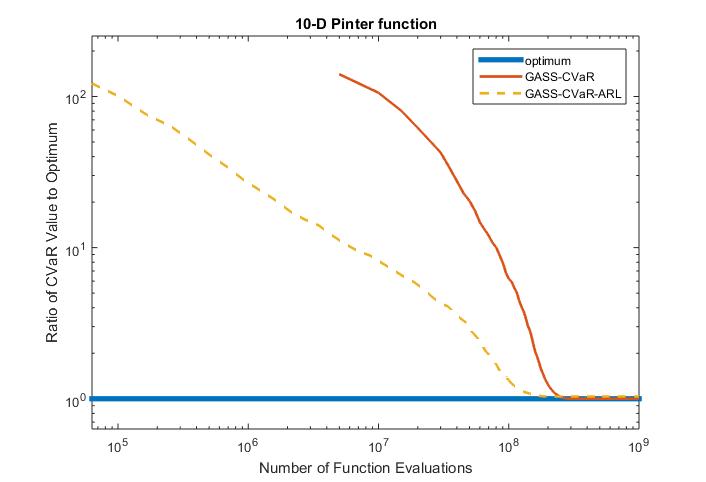}&
\hspace{-5mm}
\includegraphics[width=0.50\textwidth]{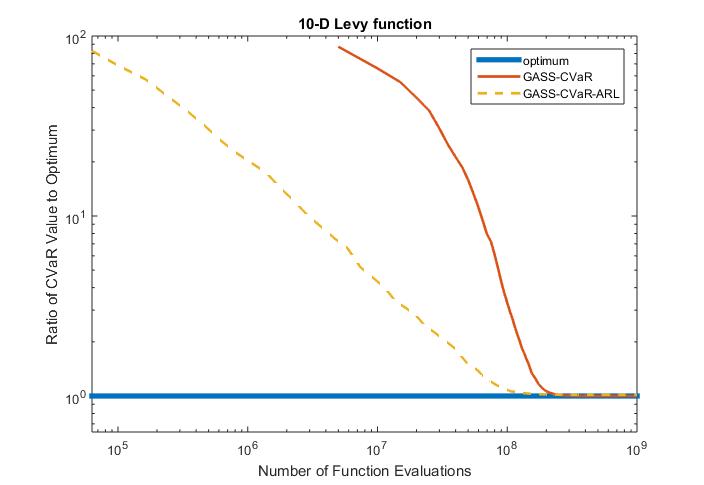}
\end{tabular}
\caption{Average Performances of GASS-CVaR and GASS-CVaR-ARL.}\label{fig.4.1}
}
\end{figure}

We run both algorithms $50$ times independently and summarize their average performance in Figure \ref{fig.4.1}.
Recall that, except for the loss function $l_0$, the minimum CVaR value is not readily available for any other loss function.
So we implement GASS-CVaR with large sample sizes $N=5\times10^3$ and $M=10^5$ to find close approximations of the true minimum CVaR values, which will be served as benchmark values later when comparing algorithm performance. Later the approximate minimum CVaR values will be used as the true CVaR values in comparing the algorithm performance.
In the upper-left plot of Figure \ref{fig.4.1} for the loss function $l_0$, the $y$-axis represents the ratio of the best CVaR values found by the algorithms to the minimum CVaR value at the target risk level $\alpha^\ast$;
for all the rest of the plots, the $y$-axis represents the same ratio, except that the minimum is replaced by the approximate minimum CVaR values from implementing GASS-CVaR with sample sizes $N$ and $M$.
We observe that both algorithms (GASS-CVaR and GASS-CVaR-ARL) perform well in finding optimal solutions and minimum CVaR values. Moreover, GASS-CVaR-ARL converges faster and often reduces the total number of function evaluations needed for convergence by $2$-$4$ times, which demonstrates the advantage of using adaptive risk levels in GASS-CVaR-ARL.

\begin{figure}[htb!]
{
\centering
\begin{tabular}{cc}
\includegraphics[width=0.50\textwidth]{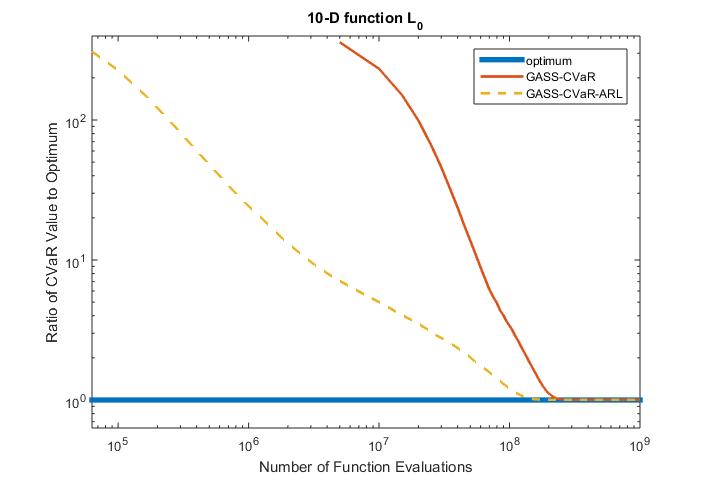}&
\hspace{-5mm}
\includegraphics[width=0.50\textwidth]{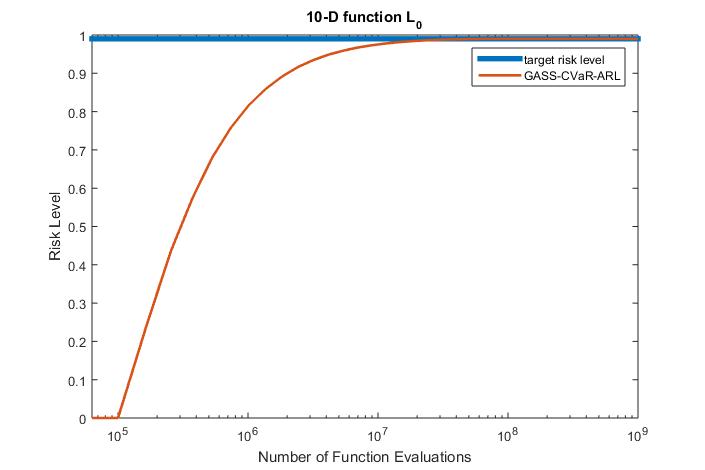}
\end{tabular}
\caption{CVaR at Mean of the Sampling Distribution and Trajectory of Risk Level.}\label{fig.4.2}
}
\end{figure}

Figure \ref{fig.4.2} includes two plots for the loss function $l_0$: the left one plots the ratio of the CVaR values evaluated at the means of the sampling distributions to the minimum CVaR value; the right one plots the trajectory of the risk level $\alpha_k$. We can see that the means of the sampling distributions in both GASS-CVaR and GASS-CVaR-ARL converge to the optimal solution, and GASS-CVaR-ARL achieves a faster convergence speed. Moreover, the risk level $\alpha_k$ in GASS-CVaR-ARL increases steadily to the target risk level $\alpha^\ast=0.99$, which indicates that the norm of the gradient decreases steadily to zero and the algorithm converges.

\section{Conclusion}\label{sec6:Conclusion}
In this paper, we study the extension of the recently proposed algorithm GASS, which is designed for the optimization of deterministic non-differentiable objectives, to the simulation optimization of risk measures such as VaR and CVaR. Instead of optimizing VaR or CVaR at the risk level of interest directly, we propose to initialize the algorithm at a small risk level, and then increase the risk level at each iteration adaptively such that the target risk level is achieved while the algorithm converges simultaneously. It enables us to adaptively reduce the number of samples needed to estimate VaR or CVaR at each iteration, leading to improvement of efficiency over the original algorithm. The numerical results demonstrate the advantage of incorporating such an adaptive updating rule on the risk level in the algorithm by showing it results in a $2$-$4$ times of total budget saving for the tested loss functions.
\section*{Acknowledgements}
This work was supported by National Science Foundation under Grants CMMI-1413790 and CAREER CMMI-1453934, and Air Force Office of Scientific Research under Grant YIP FA-9550-14-1-0059.

\appendix

\section{Proof of Theorems} \label{app:quadratic}
\begin{proof}
\textbf{Proof of Lemma \ref{lem.4.1}}.
Since $S_{\theta_k}(\cdot)$ is continuous in both ${C}_{\alpha^\ast}$ and $\gamma_{\theta_k}$, it suffices to show that for all $x\in \mathcal{X}$
\begin{equation}\label{eq.a.1}
\lim_{k\rightarrow \infty}\widehat{C}_{\alpha^\ast}(x)\rightarrow {C}_{\alpha^\ast}(x),\; w.p.1., \quad\mbox{and}\quad \lim_{k\rightarrow \infty}\widehat{\gamma}_{\theta_k}\rightarrow\gamma_{\theta_k},\; w.p.1.
\end{equation}

Let us first show the left part of the above statement. Recall that by Assumption \ref{asm.4.1}.(iii), we have $M_k\rightarrow \infty$ as $k\rightarrow \infty$. Then, we only need to show that the one-layer CVaR estimator $\widehat{C}_{\alpha^\ast}(x)$ is strongly consistent. By Lemma A.1 in \cite{zhu2016risk} this holds, where note that Assumption 3.1 in \cite{zhu2016risk} is satisfied by Assumption \ref{asm.4.2} here.

It remains to establish the right part of (\ref{eq.a.1}). In view of Assumption \ref{asm.4.1}.(ii) and \ref{asm.4.1}.(iii), we have $N_k, M_k\rightarrow\infty$ as $k\rightarrow\infty$. That is, $N_k, M_k$ go to infinity simultaneously as $k\rightarrow$. Therefore, it suffices to show
\begin{equation}\label{eq.a.2}
\lim_{N_k, M_k\rightarrow\infty} \widehat{\gamma}_{\theta_k}\rightarrow\gamma_{\theta_k},\; w.p.1.
\end{equation}
Note that
\begin{equation*}
\gamma_{\theta_k}=V_{1-\rho}(-C_{\alpha^\ast}(x)),
\end{equation*}
i.e., the $(1-\rho)$-level Value-at-Risk (VaR) of $(-C_{\alpha^\ast}(x))$ w.r.t. $f(x;\theta_k)$. Furthermore,
\begin{equation*}
\widehat{\gamma}_{\theta_k}=\widehat{V}_{1-\rho}
(-\widehat{C}_{\alpha^\ast}(x)),
\end{equation*}
i.e., the sample $(1-\rho)$-quantile of $\{-\widehat{C}_{\alpha^\ast}(x_k^i): i=1,...,N_k\}$. Therefore, $\widehat{\gamma}_{\theta_k}$ is a nested estimator of $\gamma_{\theta_k}$, where $N_k$ outer-layer samples are drawn, and for each outer-layer sample $M_k$ inner-layer samples are drawn.

Rewrite $\widehat{C}_{\alpha^\ast}(x)$ as
\begin{equation}\label{eq.a.2}
\widehat{C}_{\alpha^\ast}(x)
=C_{\alpha^\ast}(x)+\frac{1}{\sqrt{M_k}}\cdot\mathcal{E}_k(x),\quad \forall x\in\mathcal{X},
\end{equation}
where $\mathcal{E}_k(x)$ is the standardized error. By Theorem 3.3 in \cite{zhu2016risk}, we have
$$
\lim_{M_k\rightarrow\infty} \sqrt{M_k}\left(\widehat{C}_{\alpha^\ast}(x)
-C_{\alpha^\ast}(x)\right)
\overset{\mathcal{D}}\Rightarrow \mathcal{N}\left(0, \sigma^2(x)
\right),
$$
where ``$\overset{\mathcal{D}}\Rightarrow$'' denotes the convergence in distribution, and $\mathcal{N}\left(0, \sigma^2(x)\right)$ denotes a normal distribution with mean zero and variance $\sigma^2(x)$, where $\sigma^2(x)$ is the variance parameter that only depends on $x$.
Combined with (\ref{eq.a.2}), we can see that the standardized error $\mathcal{E}_k(x)$ converges to $\mathcal{N}\left(0, \sigma^2(x)\right)$ in distribution. Have establishing this, the remaining proof is identical to the proof of Theorem 3.2 in \cite{zhu2016risk}, where note that Assumption \ref{asm.4.2} here is parallel with Assumption 3.2 in \cite{zhu2016risk}.
\end{proof}

\begin{proof}\textbf{Proof of Lemma \ref{lem.4.2}}.
With a slight abuse of notation, we also use $\enorm{A}$ to denote the spectral norm of a real square matrix $A$ induced by the vector Euclidean norm. In particular, $\enorm{A}=\sqrt{\lambda_{max}(A^T A)}$, i.e., $\enorm{A}$ is the largest eigenvalue of the positive-semidefinite matrix $A^TA$. When the matrix $A$ is positive-semidefinite, $\enorm{A}$ is just the largest eigenvalue of $A$.

To facilitate the proof, let us also introduce the following notations:
\begin{equation*}
\begin{split}
&\widetilde{\mathds{Y}}_k\overset{\triangle}=\frac{1}{N_k}
\sum_{i=1}^{N_k}S_{\theta_k}\left(-C_{\alpha^\ast}(x_k^i)\right)
\Gamma(x_k^i),\;\;
\widetilde{\mathds{Z}}_k\overset{\triangle}=\frac{1}{N_k}
\sum_{i=1}^{N_k}S_{\theta_k}\left(-C_{\alpha^\ast}(x_k^i)\right),
\\
&\widehat{\mathds{Y}}_k\overset{\triangle}=\frac{1}{N_k}
\sum_{i=1}^{N_k}\widehat{S}_{\theta_k}\left(
-\widehat{C}_{\alpha^\ast}(x_k^i)\right)
\Gamma(x_k^i),\;\;
\widehat{\mathds{Z}}_k\overset{\triangle}=\frac{1}{N_k}
\sum_{i=1}^{N_k}\widehat{S}_{\theta_k}\left(
-\widehat{C}_{\alpha^\ast}(x_k^i)\right).
\end{split}
\end{equation*}
Here note that $\widetilde{\mathds{Y}}_k$, $\widehat{\mathds{Y}}_k$ are vectors because $\Gamma(\cdot)$ are vector-valued functions, and $\widetilde{\mathds{Z}}_k$, $\widehat{\mathds{Z}}_k$ are scalar-valued.

Since $C_{\alpha^\ast}(x)$ and $\Gamma(x)$ are both bounded on $\mathcal{X}$, we immediately have $|\widetilde{\mathds{Z}}_k|$ bounded below from zero and $\frac{\enorm{\widehat{\mathds{Y}}_k}}{\abs{\widehat{\mathds{Z}}_k}}
$
bounded for all $k$. Note that
\begin{eqnarray*}
b_k &=& \widehat{V}_k^{-1}\left(
\widehat{\mathbb{E}}_{q_k}[\Gamma(x)]-
\widetilde{\mathbb{E}}_{q_k}[\Gamma(x)]\right)\\
&=&\widehat{V}_k^{-1}\left(
\frac{\widehat{\mathds{Y}}_k}{\widehat{\mathds{Z}}_k}-
\frac{\widetilde{\mathds{Y}}_k}{\widetilde{\mathds{Z}}_k}\right)\\
&=&\widehat{V}_k^{-1}\left(
\frac{\widehat{\mathds{Y}}_k}{\widehat{\mathds{Z}}_k}-
\frac{\widehat{\mathds{Y}}_k}{\widetilde{\mathds{Z}}_k}+
\frac{\widehat{\mathds{Y}}_k}{\widetilde{\mathds{Z}}_k}-
\frac{\widetilde{\mathds{Y}}_k}{\widetilde{\mathds{Z}}_k}\right)\\
&=&\widehat{V}_k^{-1}\widehat{\mathds{Y}}_k\left(
\frac{\widetilde{\mathds{Z}}_k-\widehat{\mathds{Z}}_k}
{\widehat{\mathds{Z}}_k\widetilde{\mathds{Z}}_k}\right)+
\widehat{V}_k^{-1}
\frac{\widehat{\mathds{Y}}_k-\widetilde{\mathds{Y}}_k}
{\widetilde{\mathds{Z}}_k}.
\end{eqnarray*}
Therefore,
\begin{eqnarray*}
\enorm{b_k} &\le&
\frac{\enorm{\widehat{V}_k^{-1}}}{\abs{\widetilde{\mathds{Z}}_k}}
\frac{\enorm{\widehat{\mathds{Y}}_k}}{\abs{\widehat{\mathds{Z}}_k}}
\abs{\widetilde{\mathds{Z}}_k-\widehat{\mathds{Z}}_k}+
\frac{\enorm{\widehat{V}_k^{-1}}}{\abs{\widetilde{\mathds{Z}}_k}}
\abs{\widehat{\mathds{Y}}_k-\widetilde{\mathds{Y}}_k}\nonumber\\
&\le&\frac{\enorm{\widehat{V}_k^{-1}}}{\abs{\widetilde{\mathds{Z}}_k}}
\frac{\enorm{\widehat{\mathds{Y}}_k}}{\abs{\widehat{\mathds{Z}}_k}}
\frac{1}{N_k}\sum_{i=1}^{N_k}
\abs{S_{\theta_k}\left(-C_{\alpha^\ast}(x_k^i)\right)-
\widehat{S}_{\theta_k}\left(-\widehat{C}_{\alpha^\ast}(x_k^i)
\right)}\nonumber\\
&&+\frac{\enorm{\widehat{V}_k^{-1}}}{\abs{\widetilde{\mathds{Z}}_k}}
\frac{1}{N_k}\sum_{i=1}^{N_k}
\abs{\widehat{S}_{\theta_k}\left(-\widehat{C}_{\alpha^\ast}(x_k^i)
\right)-S_{\theta_k}\left(-C_{\alpha^\ast}(x_k^i)\right)}
\enorm{\Gamma(x_k^i)}.
\end{eqnarray*}

Recall that
$\widehat{V}_k=\left(
\widehat{Var}_{\theta_k}[\Gamma(x)]+\epsilon I\right)$.
Thus, it is a positive-definite matrix and its minimum eigenvalue is at least $\epsilon$. It follows that the maximum eigenvalue of $\widehat{V}_k^{-1}$ is no greater than $\epsilon^{-1}$, i.e.,
$\enorm{\widehat{V}_k^{-1}}\le \epsilon^{-1}$. Since $|\widetilde{\mathds{Z}}_k|$ is bounded below from zero, $\frac{\enorm{\widehat{\mathds{Y}}_k}}{\abs{\widehat{\mathds{Z}}_k}}$
is bounded, and $\Gamma(x)$ is bounded on $\mathcal{X}$, Lemma \ref{lem.4.1} implies that $\enorm{b_k}\rightarrow 0$ w.p.1 as $k\rightarrow \infty$.
\end{proof}

\begin{proof} \textbf{Proof of Theorem \ref{thm.4.2}}.
Let us first show the following lemma.
\begin{lemma}\label{lem.a.1}
Suppose Assumption \ref{asm.4.1} and Assumption \ref{asm.4.2} hold. Further suppose the risk level sequence $\{\alpha_k\}$ generated by (\ref{eq.3.2}) converges to the target risk level $\alpha^\ast$ w.p.1. Then the sequence $\{\theta_k\}$ generated by (\ref{eq.4.8}) converges to a limit set of the ODE (\ref{eq.4.7}) w.p.1.
\end{lemma}
Proof of Lemma \ref{lem.a.1}. Similar to the proof of Theorem \ref{thm.4.1}, we will reformulate the updating scheme (\ref{eq.4.8}) as a noisy discretization of the constrained ODE (\ref{eq.4.7}), and show both the bias and the noise are properly bounded.
Specifically, rewrite (\ref{eq.4.8}) as
\begin{equation}\label{eq.a.4}
\theta_{k+1}=\theta_k+\beta_k\left[G(\theta_k)+
\widebar{b}_k+e_k+\widebar{p}_k\right],
\end{equation}
where $G(\theta_k)$ and $e_k$ are defined as previously,
$\widebar{b}_k\overset{\triangle}=\widehat{V}_k^{-1}\left(
\widebar{\mathbb{E}}_{q_k}[\Gamma(x)]-
\widetilde{\mathbb{E}}_{q_k}[\Gamma(x)]\right)$,
and
$\widebar{p}_k$ is the projection error term that takes the current iterate back onto the constraint set $\widetilde{\Theta}$ with minimum Euclidean norm. In view of Theorem 2 in \cite{kushner2010stochastic}, it suffices to show
\begin{equation*}
\lim_{k\rightarrow \infty}\enorm{\widebar{b}_k}= 0,\quad w.p.1.
\end{equation*}
To ease the presentation, let us denote
\begin{equation*}
\widetilde{\mathbb{E}}^{\alpha}_{q_k}[\Gamma(x)]\overset{\triangle}=
\sum_{i=1}^{N_k}\frac{S_{\theta_k}
    \left(-C_{\alpha}(x_k^i)\right)}
    {\sum_{j=1}^{N_k}S_{\theta_k}
    \left(-C_{\alpha}(x_k^j)\right)}\Gamma(x^i_k).
\end{equation*}
It immediately implies that $\widetilde{\mathbb{E}}_{q_k}[\Gamma(x)]=
\widetilde{\mathbb{E}}^{\alpha^\ast}_{q_k}[\Gamma(x)]$.
Furthermore,
\begin{eqnarray}
\enorm{\widebar{b}_k}&=&\enorm{\widehat{V}_k^{-1}\left(
\widebar{\mathbb{E}}_{q_k}[\Gamma(x)]-
\widetilde{\mathbb{E}}_{q_k}[\Gamma(x)]\right)}\nonumber\\
&=&\enorm{\widehat{V}_k^{-1}\left(
\widebar{\mathbb{E}}_{q_k}[\Gamma(x)]-
\widetilde{\mathbb{E}}^{\alpha_k}_{q_k}[\Gamma(x)]\right)+
\widehat{V}_k^{-1}\left(
\widetilde{\mathbb{E}}^{\alpha_k}_{q_k}[\Gamma(x)]-
\widetilde{\mathbb{E}}^{\alpha^\ast}_{q_k}[\Gamma(x)]\right)}\nonumber\\
&\le& \enorm{\widehat{V}_k^{-1}\left(
\widebar{\mathbb{E}}_{q_k}[\Gamma(x)]-
\widetilde{\mathbb{E}}^{\alpha_k}_{q_k}[\Gamma(x)]\right)}+
\enorm{\widehat{V}_k^{-1}}
\enorm{\widetilde{\mathbb{E}}^{\alpha_k}_{q_k}[\Gamma(x)]-
\widetilde{\mathbb{E}}^{\alpha^\ast}_{q_k}[\Gamma(x)]}.\label{eq.a.5}
\end{eqnarray}
Following an argument almost identical to the proof of Lemma \ref{lem.4.2}, the first term in (\ref{eq.a.5}) converges to $0$ w.p.1 as $k\rightarrow\infty$. Note that $S_{\theta_k}(\cdot)$ is a continuous function and $C_{\alpha_k}(x)$ is continuous in $\alpha_k$. Thus, $\widetilde{\mathbb{E}}^{\alpha_k}_{q_k}[\Gamma(x)]$ is a continuous function in $\alpha_k$. Therefore, the second term in (\ref{eq.a.5}) converges to $0$ w.p.1 as $k\rightarrow\infty$ since $\enorm{\widehat{V}_k^{-1}}$ is bounded and $\alpha_k$ converges to $\alpha^\ast$ as $k\rightarrow\infty$. Proof of Lemma \ref{lem.a.1} is now complete.

In view of Lemma \ref{lem.a.1}, it remains to show that the risk level sequence $\{\alpha_k\}$ generated by (\ref{eq.3.2}) converges to the target risk level $\alpha^\ast$ w.p.1. Proof by contradiction. Since the sequence $\{\alpha_k\}$ is non-decreasing and bounded above by $\alpha^\ast$, let us assume $\lim_{k\rightarrow\infty}\alpha_k=\widebar{\alpha}^\ast$ and $\widebar{\alpha}^\ast< \alpha^\ast$ w.p.1. Conditioning on this, Lemma \ref{lem.a.1} still holds when the target risk level $\alpha^\ast$ is replaced by $\widebar{\alpha}^\ast$. That is, the algorithm GASS-CVaR-ARL converges, and the gradient sequence $\{g_k\}$ converges to $0$ w.p.1. as $k\rightarrow\infty$. Note that $g_k$ is bounded (since $\Gamma(x)$ is bounded), by bounded convergence theorem we have
\begin{equation}\label{eq.a.6}
\lim_{k\rightarrow\infty}\mathbb{E}\left[\enorm{g_k}\right]=0.
\end{equation}
Furthermore, note that
\begin{eqnarray}
\mathbb{E}\left[\enorm{\widebar{g}_k-g_k}\right]&=&
\mathbb{E}\left[\enorm{\widebar{\mathbb{E}}_{q_k}[\Gamma(x)]-
\mathbb{E}^{\widebar{\alpha}^\ast}_{q_k}[\Gamma(x)]}\right]\nonumber\\
&\le& \mathbb{E}\left[\enorm{\widebar{\mathbb{E}}_{q_k}[\Gamma(x)]-
\widetilde{\mathbb{E}}^{\widebar{\alpha}^\ast}_{q_k}[\Gamma(x)]}\right]
+\mathbb{E}\left[\enorm{\widetilde{\mathbb{E}}^
{\widebar{\alpha}^\ast}_{q_k}[\Gamma(x)]-
\mathbb{E}^{\widebar{\alpha}^\ast}_{q_k}[\Gamma(x)]}\right],
\label{eq.a.7}
\end{eqnarray}
where
\begin{equation*}
\mathbb{E}^{\widebar{\alpha}^\ast}_{q_k}[\Gamma(x)]\overset{\triangle}=
\frac{\int S_{\theta_k}\left(-C_{\widebar{\alpha}^\ast}(x)\right)
\Gamma(x)f(x;\theta_k)dx}
{\int S_{\theta_k}\left(-C_{\widebar{\alpha}^\ast}(x)\right)
f(x;\theta_k)dx}.
\end{equation*}
We have shown in the proof of Lemma \ref{lem.a.1} that
\begin{equation*}
\lim_{k\rightarrow\infty}\enorm{\widebar{\mathbb{E}}_{q_k}[\Gamma(x)]-
\widetilde{\mathbb{E}}^{\widebar{\alpha}^\ast}_{q_k}[\Gamma(x)]}=0, \quad w.p.1.
\end{equation*}
Since $\enorm{\widebar{\mathbb{E}}_{q_k}[\Gamma(x)]-
\widetilde{\mathbb{E}}^{\widebar{\alpha}^\ast}_{q_k}[\Gamma(x)]}$ is bounded, again by bounded convergence theorem
\begin{equation}\label{eq.a.8}
\lim_{k\rightarrow\infty}
\mathbb{E}\left[\enorm{\widebar{\mathbb{E}}_{q_k}[\Gamma(x)]-
\widetilde{\mathbb{E}}^{\widebar{\alpha}^\ast}_{q_k}
[\Gamma(x)]}\right]=0.
\end{equation}
Moreover, notice that $\widetilde{\mathbb{E}}^{\widebar{\alpha}^\ast}_{q_k}[\Gamma(x)]$ is a self-normalized importance sampling estimator of $\mathbb{E}^{\widebar{\alpha}^\ast}_{q_k}[\Gamma(x)]$. Applying Theorem 9.1.10 (pp. 294) in \cite{cappe2005}, we have
\begin{equation*}
\mathbb{E}\left[\abs{\widetilde{\mathbb{E}}^
{\widebar{\alpha}^\ast}_{q_k}[\Gamma_j(x)]-
\mathbb{E}^{\widebar{\alpha}^\ast}_{q_k}[\Gamma_j(x)]}^2\right]\le \frac{c_j}{N_k}, \; j=1,...,d_\theta,
\end{equation*}
where $\Gamma_j(x)$ is the $j^{th}$ element in the vector $\Gamma(x)$, and $c_j$'s are positive constants that depend on the bounds of $\Gamma_j(x)$'s on $\mathcal{X}$. Therefore, by Cauchy-Schwarz Inequality we have
\begin{equation*}
\mathbb{E}\left[
\enorm{\widetilde{\mathbb{E}}^{\widebar{\alpha}^\ast}_{q_k}[\Gamma(x)]-
\mathbb{E}^{\widebar{\alpha}^\ast}_{q_k}[\Gamma(x)]}\right]\le \sqrt{\frac{d\cdot\max_j c_j}{N_k}}.
\end{equation*}
That is,
\begin{equation}\label{eq.a.9}
\lim_{k\rightarrow\infty}
\mathbb{E}\left[
\enorm{\widetilde{\mathbb{E}}^{\widebar{\alpha}^\ast}_{q_k}[\Gamma(x)]-
\mathbb{E}^{\widebar{\alpha}^\ast}_{q_k}[\Gamma(x)]}\right]=0.
\end{equation}
Combining (\ref{eq.a.7}), (\ref{eq.a.8}) with (\ref{eq.a.9}), we have
\begin{equation*}
\lim_{k\rightarrow\infty}
\mathbb{E}\left[\enorm{\widebar{g}_k-g_k}\right]=0.
\end{equation*}
In view of (\ref{eq.a.6}), we have
\begin{equation}\label{eq.a.10}
\lim_{k\rightarrow\infty}
\mathbb{E}\left[\enorm{\widebar{g}_k}\right]=0.
\end{equation}
Since $\widebar{\alpha}^\ast<\alpha^\ast$, the sequence $\{\enorm{\widebar{g}_k}\}$ generated by (\ref{eq.3.2}) will always be above a certain positive value w.p.1 (otherwise $\alpha_k$ will converge to $\alpha^\ast$), which contradicts with (\ref{eq.a.10}). Proof is complete.
\end{proof}

\bibliographystyle{ormsv080}
\bibliography{Zhou-Bibtex}
\end{document}